\pdfoutput=1
\documentclass[11pt]{extarticle}

\usepackage[left=1in,right=1in,top=1in,bottom=1.3in]{geometry}
\usepackage{amsmath, amssymb, amsthm, latexsym, indentfirst, amsfonts, xcolor, mathtools, manfnt, subcaption, listings, longtable, textcomp, gensymb, floatpag, booktabs} 
\usepackage[shortlabels]{enumitem}
\usepackage[utf8]{inputenc} 
\usepackage[english]{babel}

\usepackage[breaklinks,backref=page]{hyperref}
\usepackage{cleveref}
\hypersetup{
	colorlinks = true, 
	urlcolor = cyan, 
	linkcolor = teal, 
	citecolor = cyan 
}

\usepackage{tikz} 
\usetikzlibrary{shapes.geometric}

\renewcommand*{\backref}[1]{}
\renewcommand*{\backrefalt}[4]{
	\ifcase #1 \textcolor{red}{Not cited.}%
	\or $\uparrow$#2%
	\else $\uparrow$#2%
	\fi%
}
\let\OLDthebibliography\thebibliography
\renewcommand\thebibliography[1]{
	\OLDthebibliography{#1}
	\setlength{\parskip}{0pt}
	\setlength{\itemsep}{0pt plus 0.3ex}
}

\newtheorem{Theorem}{Theorem}

\newtheorem{Conjecture}{Conjecture}
\newtheorem{Corollary}{Corollary}

\newtheorem{Lemma}{Lemma}

\newtheorem{Proposition}{Proposition}

\crefname{Theorem}{theorem}{theorems}
\Crefname{Theorem}{Theorem}{Theorems}
\crefname{Claim}{claim}{claims}
\Crefname{Claim}{Claim}{Claims}
\crefname{Conjecture}{conjecture}{conjectures}
\Crefname{Conjecture}{Conjecture}{Conjectures}
\crefname{Corollary}{corollary}{corollaries}
\Crefname{Corollary}{Corollary}{Corollaries}
\crefname{Example}{example}{examples}
\Crefname{Example}{Example}{Examples}
\crefname{Lemma}{lemma}{lemmas}
\Crefname{Lemma}{Lemma}{Lemmas}
\crefname{Problem}{problem}{problems}
\Crefname{Problem}{Problem}{Problems}
\crefname{Question}{question}{questions}
\Crefname{Question}{Question}{Questions}
\crefname{Proposition}{proposition}{propositions}
\Crefname{Proposition}{Proposition}{Propositions}
\crefname{Remark}{remark}{remarks}
\Crefname{Remark}{Remark}{Remarks}

\newcommand{\R}{\mathbb R}
\newcommand{\Z}{\mathbb Z}
\newcommand{\N}{\mathbb N}

\DeclareSymbolFont{endpointsymbols}{U}{txsyc}{m}{n}
\DeclareMathSymbol{\segmentdots}{\mathrel}{endpointsymbols}{"96}
\newcommand{\T}{\mathbb{T}}
\newcommand{\TrO}{\operatorname{T}}
\newcommand{\Tr}[3]{\operatorname{T}(#1,#2,#3)}
\newcommand{\TrI}[1]{\operatorname{T}(#1)}
\newcommand{\newT}{\Tr{\sqrt{3}}{\sqrt{7}}{\sqrt{7}}}
\newcommand{\rArrows}[1]{\xrightarrow{#1}}
\newcommand{\MRarrows}{\smash{\xrightarrow%
  [\smash{\raisebox{3pt}{\text{{\tiny \upshape MR}}}}]%
  {}}}

\newcommand{\trip}[2]{\left(#1,#2\right)_{\T}}
\newcommand{\lineopacity}{0.2}

\definecolor{color1}{HTML}{F09EA7}
\definecolor{color2}{HTML}{F6CA94}
\definecolor{color3}{HTML}{FAFABE}
\definecolor{color4}{HTML}{C1EBC0}
\definecolor{color5}{HTML}{C7CAFF}
\definecolor{color6}{HTML}{CDABEB}
\definecolor{color7}{HTML}{F6C2F3}
\colorlet{strongred}{red!90!black}
\colorlet{strongblue}{blue!90!black}

\tikzset{
  hexagon/.style={
    draw,
    regular polygon,
    regular polygon sides=6,
    shape border rotate=30,
    minimum size={2*\hexsize cm},
    inner sep=0pt,
    font=\small,
    align=center
  }
}

\newcommand{\tricoord}[2]{
  {\trisize * (#1 + 0.5* (#2 + 1))},
   {(sqrt(3)/2*\trisize*(#2 + 1))}
}

\tikzstyle{vector}=[->,very thick]
\usetikzlibrary{arrows.meta, calc}

\newcommand{\trigrid}[2]{
  \def\nc{#2}
  \def\nr{#1}
  \pgfmathsetmacro\s{\trisize}
        
        \pgfmathsetmacro\nx{\nc + 1}
        \pgfmathsetmacro\ny{ceil(\nr * 0.5)}
        \pgfmathsetmacro\nyy{(2+\nr)*sin(60)*\s}
        \clip (0.2*\s,0.6*\s) rectangle (\s*\nx -\s*0.2,\nyy -\s*0.6);

        \foreach \j in {0,...,\ny} {
            \foreach \i in {0,...,\nx} {
               \draw[dotted] (0:\i *\s)++(90:{\s *(1+2*\j)*sin(60)})--++(\s,0);
               \draw[dotted] (0:\s *\i) ++(60:\s *\j) ++(120:\s*\j) {} --++(60:\s *2) {} --++(-1*\s,0) {} --++(-60:\s) {} --++(-60:\s) {};
            }}}

\usepackage{pgfplots}
\pgfplotsset{compat=1.18}

\usepgfplotslibrary{external,fillbetween}

\pgfplotsset{
    myplot/.style={
        axis lines=middle,
        axis line style={black, thick},
        enlargelimits=false,
        xmin=-0.05, xmax=1.05,
        ymin=-0.05, ymax=1.05,
        xtick={0,1},
        ytick={0,1},
        xlabel={},
        ylabel={},
        width=10cm,
        height=10cm,
        axis equal image,
        clip=true,
    }
}

\definecolor{codegreen}{rgb}{0,0.6,0}
\definecolor{codegray}{rgb}{0.5,0.5,0.5}
\definecolor{codepurple}{rgb}{0.58,0,0.82}
\definecolor{backcolour}{rgb}{0.95,0.95,0.92}

\lstdefinestyle{mystyle}{
    backgroundcolor=\color{backcolour},   
    commentstyle=\color{codegreen},
    keywordstyle=\color{magenta},
    numberstyle=\tiny\color{codegray},
    stringstyle=\color{codepurple},
    basicstyle=\ttfamily\footnotesize,
    breakatwhitespace=false,         
    breaklines=true,                 
    captionpos=b,                    
    keepspaces=true,                 
    numbers=left,                    
    numbersep=5pt,                  
    showspaces=false,                
    showstringspaces=false,
    showtabs=false,                  
    tabsize=2
}

\title{Monochromatic or rainbow triangles in the triangular lattice\footnote{This paper is a shortened version of the first author's master's thesis completed in July 2026.}}
\author{Till Rauch-M\"uller\thanks{Karlsruhe Institute of Technology, Karlsruhe, Germany, \href{mailto:till.rauch@me.com}{\tt till.rauch@me.com}.}
\and Arsenii Sagdeev\thanks{Karlsruhe Institute of Technology, Karlsruhe, Germany, \href{mailto:sagdeevarsenii@gmail.com}{\tt sagdeevarsenii@gmail.com}.}}
\date{}

\begin{document}

\maketitle

\begin{abstract}
   Erd\H{o}s, Graham, Montgomery, Rothschild, Spencer, and Straus conjectured that for each non-equilateral triangle $\operatorname{T}$, every 2-coloring of the plane contains a congruent copy of $\operatorname{T}$ with monochromatic vertex set. We confirm this conjecture for triangles with a ratio $\sqrt{7/3}$ between two side lengths. We further prove that every coloring of the plane, regardless of the number of colors used, contains congruent copies of 36 specific triangles with monochromatic or rainbow vertex sets, partially answering a question of Geh\'er, Sagdeev, and T\'oth. Both these results were obtained by considering finite configurations on the triangular lattice. 
\end{abstract}

\noindent
\textbf{AI disclosure.} No AI assistance was used in obtaining the results of this paper and in its writing.

\section{Introduction}

For a set $S \subset \R^2$ and a positive integer $r$, we write $\R^2 \xrightarrow{r} S$ if every $r$-coloring of $\R^2$ contains a monochromatic congruent copy of $S$; otherwise, we write $\R^2 \not \xrightarrow{r} S$. The classical Hadwiger--Nelson problem asks for the maximum $r$ such that $\R^2 \xrightarrow{r} \ \segmentdots$, where $\segmentdots$ denotes a set of two points at distance 1, see Soifer's book~\cite{soifer}  for a history of this problem. In 2018, de Grey~\cite{degrey} proved that $\R^2 \xrightarrow{4} \ \segmentdots$, see also Exoo and Ismailescu~\cite{exoo2020chromatic} and Parts~\cite{parts}; all these proofs are computer-assisted. (In a recent, not peer-reviewed manuscript~\cite{AI26}, OpenAI claims that $\R^2 \xrightarrow{5} \ \segmentdots$.) From the other direction, Hadwiger's classical `hexagonal' coloring~\cite{hadwiger} shows that $\R^2 \not\xrightarrow{7} \segmentdots$.

In the present paper, instead of the set $\segmentdots$, we consider \textit{triangles}, i.e., 3-element sets. Let $\Tr{a}{b}{c}$ be a set of three points with pairwise distances $a,b$, and $c$; in what follows, we assume that $a,b$, and $c$ are positive reals satisfying the triangular inequality but they are not necessarily monotonically ordered. Erd\H{o}s, Graham, Montgomery, Rothschild, Spencer, and Straus initiated this line of research in~\cite{EGMRSS1, EGMRSS2, EGMRSS3} and made the following fundamental conjecture.

\begin{Conjecture}[\cite{EGMRSS3}]
    \label{con:chromLB}
    If $\TrO$ is a non-equilateral triangle, then $\R^2 \xrightarrow{2} \TrO$.
\end{Conjecture}

Despite more than 50 years of research, \Cref{con:chromLB} remains confirmed only for a few special classes of triangles, see \Cref{sec:mono} for a somewhat technical description. In the present paper, we find a new such class.

\begin{Theorem} \label{th:1}
    If $\frac{b}{a}=\sqrt{\frac{7}{3}}$, then $\R^2 \xrightarrow{2} \Tr{a}{b}{c}$.
\end{Theorem}

In addition, using a computational approach, we confirm \Cref{con:chromLB} for 32 new specific triangles, see Sections~\ref{sec:mono} and~\ref{sec:comp} for the context and illustrative figures. Let

\vspace{-5mm}
\begin{alignat*}{4}
\mathcal{A}_{\text{\upshape M}}=\bigl\{
&\Tr{\sqrt{3}}{\sqrt{7}}{\sqrt{7}},\; &&\Tr{1}{3}{\sqrt{13}},\; &&\Tr{\sqrt{7}}{2 \sqrt{3}}{\sqrt{13}},\; &&\Tr{\sqrt{3}}{\sqrt{7}}{\sqrt{19}},\; \\
&\Tr{1}{\sqrt{13}}{\sqrt{21}},\; &&\Tr{2}{\sqrt{19}}{\sqrt{21}},\; &&\Tr{3}{\sqrt{13}}{2 \sqrt{7}},\; &&\Tr{\sqrt{3}}{\sqrt{19}}{2 \sqrt{7}},\; \\
&\Tr{1}{\sqrt{21}}{2 \sqrt{7}},\; &&\Tr{2 \sqrt{3}}{\sqrt{13}}{\sqrt{31}},\; &&\Tr{2}{\sqrt{21}}{\sqrt{31}},\; &&\Tr{3}{2 \sqrt{7}}{\sqrt{31}},\; \\
&\Tr{\sqrt{7}}{\sqrt{13}}{6},\; &&\Tr{\sqrt{7}}{\sqrt{19}}{6},\; &&\Tr{\sqrt{3}}{\sqrt{19}}{\sqrt{37}},\; &&\Tr{2}{\sqrt{21}}{\sqrt{37}},\; \\
&\Tr{\sqrt{13}}{\sqrt{21}}{\sqrt{37}},\; &&\Tr{\sqrt{3}}{2 \sqrt{7}}{\sqrt{37}},\; &&\Tr{2 \sqrt{3}}{\sqrt{31}}{\sqrt{37}},\; &&\Tr{\sqrt{7}}{6}{\sqrt{37}},\; \\
&\Tr{\sqrt{7}}{\sqrt{19}}{\sqrt{39}},\; &&\Tr{1}{2 \sqrt{7}}{\sqrt{39}},\; &&\Tr{2 \sqrt{3}}{\sqrt{13}}{\sqrt{43}},\; &&\Tr{\sqrt{7}}{\sqrt{19}}{4 \sqrt{3}},\; \\
&\Tr{\sqrt{7}}{\sqrt{31}}{4 \sqrt{3}},\; &&\Tr{\sqrt{7}}{\sqrt{43}}{4 \sqrt{3}},\; &&\Tr{2 \sqrt{3}}{\sqrt{13}}{7},\; &&\Tr{\sqrt{13}}{\sqrt{21}}{7},\; \\
&\Tr{\sqrt{3}}{2 \sqrt{7}}{7},\; &&\Tr{\sqrt{3}}{\sqrt{37}}{7},\; &&\Tr{\sqrt{3}}{\sqrt{37}}{2 \sqrt{13}},\; &&\Tr{\sqrt{7}}{\sqrt{37}}{3 \sqrt{7}}
\bigr\}.
\end{alignat*}

\begin{Theorem} \label{th:2}
    If $\TrO \in \mathcal{A}_{\text{\upshape M}}$, then $\R^2 \xrightarrow{2} \TrO$.
\end{Theorem}

From the other direction, Graham~\cite{graham} conjectured that $\R^2\not\xrightarrow{3} \TrO$ for each triangle $\TrO$. Building on the work of Aichholzer and Perz~\cite{aich}, Mundinger, Zimmer, Kiem, Spiegel, and Pokutta in their computer-assisted work~\cite{mund} confirmed this conjecture for a large family of triangles via explicit $3$-colorings. In general, however, no better result than the trivial $\R^2\not\xrightarrow{7} \TrO$ is known.

Next, we consider a `canonical' variant of these problems.  For a set $S \subset \R^n$, we write $\R^n \MRarrows S$ if every coloring of $\R^n$, regardless of the number of colors used, contains a congruent copy of $S$ that is either monochromatic or \textit{rainbow}, i.e., with points of pairwise different colors. This  variant was introduced by Mao, Ozeki, and Wang~\cite{mao}, who proved that $\R^3 \MRarrows \Tr{1}{\sqrt{3}}{2}$. Cheng and Xu~\cite{cheng2025euclidean} as well as Geh\'er, Sagdeev, and T\'oth~\cite{geher} strengthened this result by showing that $\R^3 \MRarrows \TrO$ if $\TrO$ is a right or an acute triangle, respectively. Fang, Ge, Shu, Xu, Xu, and Yang~\cite{fang} proved that if $\TrO$ is an obtuse triangle, then $\R^4 \MRarrows \TrO$. Geh\'er, Sagdeev, and T\'oth explicitly asked in~\cite{geher} if the dimension in these results can be further reduced to 2 for non-equilateral triangles. It is somewhat natural to expect such a strengthening of \Cref{con:chromLB} since it seems hard to avoid rainbow copies of a given non-equilateral triangle in the plane if more than 2 colors are used.

\begin{Conjecture}
    \label{con:MR}
    If $\TrO$ is a non-equilateral triangle, then $\R^2 \MRarrows \TrO$.
\end{Conjecture}

In the present paper, using a computational approach, we take a first step towards \Cref{con:MR} by confirming it for 36 specific triangles. Let 

\vspace{-5mm}
\begin{alignat*}{4}
\mathcal{A}_{\text{\upshape MR}}=\bigl\{
&\Tr{1}{1}{\sqrt{3}},\; &&\Tr{1}{\sqrt{3}}{2},\; &&\Tr{1}{\sqrt{3}}{\sqrt{7}},\; &&\Tr{\sqrt{3}}{2}{\sqrt{7}},\; \\
&\Tr{1}{\sqrt{7}}{3},\; &&\Tr{1}{\sqrt{7}}{2 \sqrt{3}},\; &&\Tr{\sqrt{3}}{2}{\sqrt{13}},\; &&\Tr{\sqrt{3}}{\sqrt{7}}{\sqrt{13}},\; \\
&\Tr{1}{3}{\sqrt{13}},\; &&\Tr{1}{2 \sqrt{3}}{\sqrt{13}},\; &&\Tr{\sqrt{7}}{2 \sqrt{3}}{\sqrt{13}},\; &&\Tr{\sqrt{3}}{\sqrt{7}}{4},\; \\
&\Tr{\sqrt{7}}{3}{\sqrt{19}},\; &&\Tr{1}{2 \sqrt{3}}{\sqrt{19}},\; &&\Tr{\sqrt{3}}{\sqrt{13}}{\sqrt{19}},\; &&\Tr{2}{\sqrt{7}}{\sqrt{21}},\; \\
&\Tr{1}{\sqrt{13}}{\sqrt{21}},\; &&\Tr{2}{\sqrt{13}}{\sqrt{21}},\; &&\Tr{1}{\sqrt{19}}{\sqrt{21}},\; &&\Tr{2}{\sqrt{19}}{\sqrt{21}},\; \\
&\Tr{3}{\sqrt{13}}{2 \sqrt{7}},\; &&\Tr{\sqrt{3}}{\sqrt{19}}{2 \sqrt{7}},\; &&\Tr{\sqrt{7}}{2 \sqrt{3}}{\sqrt{31}},\; &&\Tr{2 \sqrt{3}}{\sqrt{13}}{\sqrt{31}},\; \\
&\Tr{2}{\sqrt{21}}{\sqrt{31}},\; &&\Tr{3}{2 \sqrt{7}}{\sqrt{31}},\; &&\Tr{\sqrt{3}}{\sqrt{19}}{\sqrt{37}},\; &&\Tr{2}{\sqrt{21}}{\sqrt{37}},\; \\
&\Tr{\sqrt{3}}{\sqrt{31}}{\sqrt{37}},\; &&\Tr{2 \sqrt{3}}{\sqrt{31}}{\sqrt{37}},\; &&\Tr{1}{2 \sqrt{7}}{\sqrt{39}},\; &&\Tr{2 \sqrt{3}}{\sqrt{19}}{\sqrt{43}},\; \\
&\Tr{1}{\sqrt{37}}{4 \sqrt{3}},\; &&\Tr{\sqrt{3}}{\sqrt{37}}{7},\; &&\Tr{\sqrt{3}}{\sqrt{37}}{2 \sqrt{13}},\; &&\Tr{1}{4 \sqrt{3}}{\sqrt{61}}
\bigr\}.
\end{alignat*}

\enlargethispage{4mm}
\begin{Theorem} \label{th:3}
    If $\TrO \in \mathcal{A}_{\text{\upshape MR}}$, then $\R^2 \MRarrows \TrO$.
\end{Theorem}
\newpage

Note that the set $\mathcal{A}_{\text{\upshape MR}}$ contains 29 obtuse, 3 right, and 4 acute triangles. In \Cref{sec:5}, we also provide a self-contained proof for the right triangle $\Tr{1}{\sqrt{3}}{2}\in \mathcal{A}_{\text{\upshape MR}}$, considered by Mao, Ozeki, and Wang~\cite{mao}, to illustrate the additional difficulty arising when more than 2 colors are allowed.

Finally, for recent developments in Euclidean Ramsey theory, we refer the reader to the papers~\cite{ALS, conlon2026non, currier2026avoiding, davies2026measurable, davies2026prime, didin2026chromatic, ge2026all, kovavc2026coloring, kupavskii2026cutting,  moore, palvolgyi2026nearcircumsphere, palvolgyi2026cyclic, shaw2026regular, sokolov2026chromatic, tomon2026geometric} and the references therein.

\vspace{2mm}
\noindent
\textbf{Paper outline.} In \Cref{sec:mono}, we survey previous results towards Conjecture~\ref{con:chromLB}. \Cref{sec:comp} contains computer-assisted proofs of Theorems~\ref{th:2} and~\ref{th:3}, see Appendices~\ref{sec:appA} and~\ref{sec:app} for the supplementary data and the program code, respectively. In \Cref{sec:4}, we prove \Cref{th:1} without computer assistance. In \Cref{sec:5} we give a self-contained proof of the special case $\R^2\MRarrows \Tr{1}{\sqrt{3}}{2}$ of \Cref{th:3}. We conclude the paper with final remarks in \Cref{sec:concl}.

\vspace{2mm}
\noindent
\textbf{Note added.} After this writeup has been completed, we learned from 
Anton Bernshteyn and Hemanshu Kaul that a complete solution to \Cref{con:chromLB} was announced in a github repository~\cite{AI26_2} on the 19th of October 2026. The repository does not even contain a pdf-file with the proof, rather a sketch and a lean-verification code. We emphasize that this paper is a shortened version of the first author's master's thesis completed in July 2026. See also the footnote in \Cref{sec:concl}.

\section{Previous results towards Conjecture~\ref{con:chromLB}} \label{sec:mono}

Erd\H{o}s, Graham, Montgomery, Rothschild, Spencer, and Straus conjectured~\cite{EGMRSS3} that $\R^2 \xrightarrow{2} \TrO$ for each non-equilateral\footnote{For an equilateral triangle $\TrO$ of height $h$, they observed that a 2-coloring of the plane with parallel half-open color-alternating strips of width $h$ shows that $\R^2 \not\xrightarrow{2} \TrO$. They also conjectured that such a coloring is unique, except for some freedom on the boundaries. That was however disproved by Jelínek, Kynčl, Stolař, and Valla~\cite{Jelinek}.} triangle $\TrO$ and proved this for several classes of triangles.

\begin{Theorem}[\cite{EGMRSS3}]
  If $\TrO=\Tr{a}{b}{c}$ satisfies one of the following conditions, then $\R^2\rArrows{2}\TrO$: 
  \begin{enumerate}
  \item[\normalfont (i)]{$\TrO$ has an angle of $30\degree$ or, equivalently, $c^2=a^2+b^2-\sqrt{3}ab$,}
  \item[\normalfont (ii)]{$\TrO$ has an angle of $150\degree$ or, equivalently, $c^2=a^2+b^2+\sqrt{3}ab$,}
  \item[\normalfont (iii)]{$\TrO$ has angles $(\alpha, 2\alpha, 180\degree -3\alpha)$ for $0 < \alpha < 60\degree$ or, equivalently, $b^2=a^2+ac$,\,\footnote{Theorem~9 in~\cite{EGMRSS3} uses an additional assumption $\alpha \neq 45\degree$, the case $\alpha = 45\degree$ is covered separately by \cite[Theorem~14]{EGMRSS3}. This class also contains degenerate triangles of the form $\Tr{a}{2a}{3a}$ corresponding to the limit $\alpha\to 0$.}}
  \item[\normalfont (iv)]{$\TrO$ has angles $(180\degree - \alpha, 180\degree - 2\alpha, 3\alpha - 180\degree)$ for $60\degree < \alpha < 90\degree$, or, equivalently, $b^2=a^2-ac$,}
  \item[\normalfont (v)]{\(a^6 - 2a^4b^2 + a^2b^4 - 3a^2b^2c^2 + b^2c^4 = 0,\)\,\footnote{In~\cite[Introduction]{EGMRSS3}, the last term is incorrectly given as $b^2c^2$, which may explain subsequent errors in \cite{graham}.}}
  \item[\normalfont (vi)]{
    \(a^4c^2 + a^2b^4 - 5a^2b^2c^2 + b^2c^4 = 0,\)}
  \item[\normalfont (vii)]{\(4a^2c^2 -b^2c^2 - a^4 = 0,\)}
  \item[\normalfont (viii)]{\(c^2 = 2a^2 + b^2.\)\,\footnote{This result is mentioned in Section~3 of~\cite{EGMRSS3} but not in the introduction, which may explain its rediscovery by Shader~\cite{shader} and its attribution to Shader by Graham~\cite{graham}.}}
  \end{enumerate}
\end{Theorem}

\newpage
In 1976, Shader~\cite{shader} proved \Cref{con:chromLB} for two new classes of triangles.

\begin{Theorem}[\cite{shader}]
  If $\TrO=\Tr{a}{b}{c}$ satisfies one of the following conditions, then $\R^2\rArrows{2}\TrO$: 
  \begin{enumerate}
  \item[\normalfont (ix)]{$\TrO$ has an angle of $90\degree$ or, equivalently, $c^2=a^2+b^2$,}
  \item[\normalfont (x)]{\(4b^2=a^2 + c^2.\)\,\footnote{Shader~\cite{shader} assumed an additional condition $\sqrt{3/2} < a/b < \sqrt{5/2}$, though the proof applies verbatim even without it, as Yan observed in a bachelor's thesis~\cite{yan}.}}
  \end{enumerate}
\end{Theorem}

All known proofs of \Cref{con:chromLB} for infinite classes of triangles rely on the following fundamental result from~\cite{EGMRSS3} or on its corollary for isosceles triangles. 
We abbreviate $\Tr{a}{a}{a}$ by $\TrO(a)$.

\begin{Theorem}[\cite{EGMRSS3}]
  \label{theorem:equi}
  A given $2$-coloring of the plane contains a monochromatic copy of $\Tr{a}{b}{c}$ if and only if it contains a monochromatic copy of either $\TrI{a}$, or $\TrI{b}$, or $\TrI{c}$.
\end{Theorem}

\begin{Corollary}[\cite{EGMRSS3}]
  \label{lem:ratio}
  If $\R^2\rArrows{2} \Tr{a}{a}{b}$, then $\R^2\rArrows{2} \Tr{a}{b}{c}$ for all $c$, $|a-b|\le c \le a+b$.
\end{Corollary}

For each base-to-side ratio of the isosceles triangles in classes (i)–(xii), \Cref{lem:ratio} confirms \Cref{con:chromLB} for all triangles with two side lengths in that ratio.

\begin{Corollary}[\cite{graham}]
    We have $\R^2\rArrows{2} \Tr{a}{b}{c}$ if
    \begin{enumerate}
\item[\normalfont (xi)]{$\frac{b}{a}=2\sin{\frac{\theta}{2}}$ for $\theta\in\left\{30\degree, 36\degree, 90\degree, 120\degree\right\}$\footnote{The corresponding lists in~\cite{EGMRSS3} and~\cite{graham} give the angle $\theta=36^\degree$ incorrectly as $72\degree$. In Section~3 of~\cite{EGMRSS3}, there are also two additional angles $\theta = 108^\degree $ and $\theta = 150^\degree$, omitted from the introduction; these are in fact redundant since $\frac{b}{a}$ equals $2\sin{\frac{108\degree}{2}}$ or $2\sin{\frac{150\degree}{2}}$ if and only if $\frac{a}{b}$ equals $2\sin{\frac{36\degree}{2}}$ or $2\sin{\frac{30\degree}{2}}$, respectively.}}
\end{enumerate}
\end{Corollary}

For more than four decades, there was no further (\textit{unconditional},  see Shkredov~\cite{shkredov2015some}) progress on \Cref{con:chromLB}, until in 2024, Currier, Moore, and Yip~\cite{currier} confirmed it for 3-term arithmetic progressions, i.e., degenerate triangles of the form $\Tr{a}{a}{2a}$. Combining this with \Cref{lem:ratio}, they extended the result to a new class of triangles.

\begin{Theorem}[\cite{currier}]
    We have $\R^2\rArrows{2} \Tr{a}{b}{c}$ if
  \begin{enumerate}
    \item[ \normalfont (xii)]{$\frac{b}{a}=2=2\sin{\frac{180\degree}{2}}$.}
  \end{enumerate}
\end{Theorem}

\begin{figure}[p]
  \thisfloatpagestyle{empty}
  \begin{subfigure}[b]{0.32\textwidth}
    \centering
    \resizebox{\linewidth}{!}{
      \begin{tikzpicture}
        \begin{axis}[myplot]
          \addplot[name path=f, domain=.5:1, samples=200, red, draw=none] {sqrt(1-x^2)};
          \path[name path=axis] (axis cs:0,0) -- (axis cs:1,0);
          \addplot[
            thick,
            color=blue,
            fill=blue, 
            fill opacity=0.05
          ]
          fill between[
            of=f and axis,
            soft clip={domain=.5:1}
          ];
          \addplot [
            thin,
            blue,
            mark=none
          ] table {curve_data/30deg.dat};
        \end{axis}
      \end{tikzpicture}
    }
    \caption*{(i)}
  \end{subfigure}
  \hfill
  \begin{subfigure}[b]{0.32\textwidth}
    \centering
    \resizebox{\linewidth}{!}{
      \begin{tikzpicture}
        \begin{axis}[myplot]
          \addplot[name path=f, domain=.5:1, samples=200, red, draw=none] {sqrt(1-x^2)};
          \path[name path=axis] (axis cs:0,0) -- (axis cs:1,0);
          \addplot[
            thick,
            color=blue,
            fill=blue, 
            fill opacity=0.05
          ]
          fill between[
            of=f and axis,
            soft clip={domain=.5:1}
          ];
          \addplot [
            thin,
            blue,
            mark=none
          ] table {curve_data/150deg.dat};
        \end{axis}
      \end{tikzpicture}
    }
    \caption*{(ii)}
  \end{subfigure}
  \hfill
  \begin{subfigure}[b]{0.32\textwidth}
    \centering
    \resizebox{\linewidth}{!}{
      \begin{tikzpicture}
        \begin{axis}[myplot]
          \addplot[name path=f, domain=.5:1, samples=200, red, draw=none] {sqrt(1-x^2)};
          \path[name path=axis] (axis cs:0,0) -- (axis cs:1,0);
          \addplot[
            thick,
            color=blue,
            fill=blue, 
            fill opacity=0.05
          ]
          fill between[
            of=f and axis,
            soft clip={domain=.5:1}
          ];
          \addplot [
            thin,
            blue,
            mark=none
          ] table {curve_data/a_2a.dat};
        \end{axis}
      \end{tikzpicture}
    }
    \caption*{(iii)}
  \end{subfigure}
  \par
  \begin{subfigure}[b]{0.32\textwidth}
    \centering
    \resizebox{\linewidth}{!}{
      \begin{tikzpicture}
        \begin{axis}[myplot]
          \addplot[name path=f, domain=.5:1, samples=200, red, draw=none] {sqrt(1-x^2)};
          \path[name path=axis] (axis cs:0,0) -- (axis cs:1,0);
          \addplot[
            thick,
            color=blue,
            fill=blue, 
            fill opacity=0.05
          ]
          fill between[
            of=f and axis,
            soft clip={domain=.5:1}
          ];
          \addplot [
            thin,
            blue,
            mark=none
          ] table {curve_data/180minusA.dat};
        \end{axis}
      \end{tikzpicture}
    }
    \caption*{(iv)}
  \end{subfigure}
  \hfill
  \begin{subfigure}[b]{0.32\textwidth}
    \centering
    \resizebox{\linewidth}{!}{
      \begin{tikzpicture}
        \begin{axis}[myplot]
          \addplot[name path=f, domain=.5:1, samples=200, red, draw=none] {sqrt(1-x^2)};
          \path[name path=axis] (axis cs:0,0) -- (axis cs:1,0);
          \addplot[
            thick,
            color=blue,
            fill=blue, 
            fill opacity=0.05
          ]
          fill between[
            of=f and axis,
            soft clip={domain=.5:1}
          ];
          \addplot [
            thin,
            blue,
            mark=none
          ] table {curve_data/comp1_1.dat};
          \addplot [
            thin,
            blue,
            mark=none
          ] table {curve_data/comp1_2.dat};
        \end{axis}
      \end{tikzpicture}
    }
    \caption*{(v)}
  \end{subfigure}
  \hfill
  \begin{subfigure}[b]{0.32\textwidth}
    \centering
    \resizebox{\linewidth}{!}{
      \begin{tikzpicture}
        \begin{axis}[myplot]
          \addplot[name path=f, domain=.5:1, samples=200, red, draw=none] {sqrt(1-x^2)};
          \path[name path=axis] (axis cs:0,0) -- (axis cs:1,0);
          \addplot[
            thick,
            color=blue,
            fill=blue, 
            fill opacity=0.05
          ]
          fill between[
            of=f and axis,
            soft clip={domain=.5:1}
          ];
          \addplot [
            thin,
            blue,
            mark=none
          ] table {curve_data/comp2.dat};
        \end{axis}
      \end{tikzpicture}
    }
    \caption*{(vi)}
  \end{subfigure}
  \par
  \begin{subfigure}[b]{0.32\textwidth}
    \centering
    \resizebox{\linewidth}{!}{
      \begin{tikzpicture}
        \begin{axis}[myplot]
          \addplot[name path=f, domain=.5:1, samples=200, red, draw=none] {sqrt(1-x^2)};
          \path[name path=axis] (axis cs:0,0) -- (axis cs:1,0);
          \addplot[
            thick,
            color=blue,
            fill=blue, 
            fill opacity=0.05
          ]
          fill between[
            of=f and axis,
            soft clip={domain=.5:1}
          ];
          \addplot [
            thin,
            blue,
            mark=none
          ] table {curve_data/isoc.dat};
        \end{axis}
      \end{tikzpicture}
    }
    \caption*{(vii)}
  \end{subfigure}
  \hfill
  \begin{subfigure}[b]{0.32\textwidth}
    \centering
    \resizebox{\linewidth}{!}{
      \begin{tikzpicture}
        \begin{axis}[myplot]
          \addplot[name path=f, domain=.5:1, samples=200, red, draw=none] {sqrt(1-x^2)};
          \path[name path=axis] (axis cs:0,0) -- (axis cs:1,0);
          \addplot[
            thick,
            color=blue,
            fill=blue, 
            fill opacity=0.05
          ]
          fill between[
            of=f and axis,
            soft clip={domain=.5:1}
          ];
          \addplot [
            thin,
            blue,
            mark=none
          ] table {curve_data/aSquaredPlus2BSquared.dat};
        \end{axis}
      \end{tikzpicture}
    }
    \caption*{(viii)}
  \end{subfigure}
  \hfill
  \begin{subfigure}[b]{0.32\textwidth}
    \centering
    \resizebox{\linewidth}{!}{
      \begin{tikzpicture}
        \begin{axis}[myplot]
          \addplot[name path=f, domain=.5:1, samples=200, red, draw=none] {sqrt(1-x^2)};
          \path[name path=axis] (axis cs:0,0) -- (axis cs:1,0);
          \addplot[
            thick,
            color=blue,
            fill=blue, 
            fill opacity=0.05
          ]
          fill between[
            of=f and axis,
            soft clip={domain=.5:1}
          ];
          \addplot [
            thin,
            blue,
            mark=none
          ] table {curve_data/90deg.dat};
        \end{axis}
      \end{tikzpicture}
    }
    \caption*{(ix)}
  \end{subfigure}
  \par
  \begin{subfigure}[t]{0.32\textwidth}
    \centering
    \resizebox{\linewidth}{!}{
      \begin{tikzpicture}
        \begin{axis}[myplot]
          \addplot[name path=f, domain=.5:1, samples=200, red, draw=none] {sqrt(1-x^2)};
          \path[name path=axis] (axis cs:0,0) -- (axis cs:1,0);
          \addplot[
            thick,
            color=blue,
            fill=blue, 
            fill opacity=0.05
          ]
          fill between[
            of=f and axis,
            soft clip={domain=.5:1}
          ];
          \addplot [
            thin,
            blue,
            mark=none
          ] table {curve_data/bach_diff.dat};
          \addplot [
            red,
            thick,
            dashed,
            mark=none
          ] table {curve_data/prev_bach.dat};
        \end{axis}
      \end{tikzpicture}
    }
    \captionsetup{justification=centering}
    \caption*{(x) Shader's subclass dashed, \\ Yan's extension solid.}
  \end{subfigure}
  \hfill
  \begin{subfigure}[t]{0.32\textwidth}
    \centering
    \resizebox{\linewidth}{!}{
      \begin{tikzpicture}
        \begin{axis}[myplot]
          \addplot[name path=f, domain=.5:1, samples=200, red, draw=none] {sqrt(1-x^2)};
          \path[name path=axis] (axis cs:0,0) -- (axis cs:1,0);
          \addplot[
            thick,
            color=blue,
            fill=blue, 
            fill opacity=0.05
          ]
          fill between[
            of=f and axis,
            soft clip={domain=.5:1}
          ];
          \addplot [
            thin,
            blue,
            mark=none
          ] table {curve_data/ratio_30deg.dat};
          \addplot [
            thin,
            blue,
            mark=none
          ] table {curve_data/ratio_36deg.dat};
          \addplot [
            thin,
            blue,
            mark=none
          ] table {curve_data/ratio_90deg.dat};
          \addplot [
            thin,
            blue,
            mark=none
          ] table {curve_data/ratio_120deg.dat};
        \end{axis}
      \end{tikzpicture}
    }
    \caption*{(xi)}
  \end{subfigure}
  \hfill
  \begin{subfigure}[t]{0.32\textwidth}
    \centering
    \resizebox{\linewidth}{!}{
      \begin{tikzpicture}
        \begin{axis}[myplot]
          \addplot[name path=f, domain=.5:1, samples=200, red, draw=none] {sqrt(1-x^2)};
          \path[name path=axis] (axis cs:0,0) -- (axis cs:1,0);
          \addplot[
            thick,
            color=blue,
            fill=blue, 
            fill opacity=0.05
          ]
          fill between[
            of=f and axis,
            soft clip={domain=.5:1}
          ];
          \addplot [
            thin,
            blue,
            mark=none
          ] table {curve_data/degen.dat};
        \end{axis}
      \end{tikzpicture}
    }
    \caption*{(xii)}
  \end{subfigure}
  \par
  \caption{Blue lines consist of points identifying normed triangles in classes (i)-(xii).} 
  \label{fig:classes}
\end{figure}
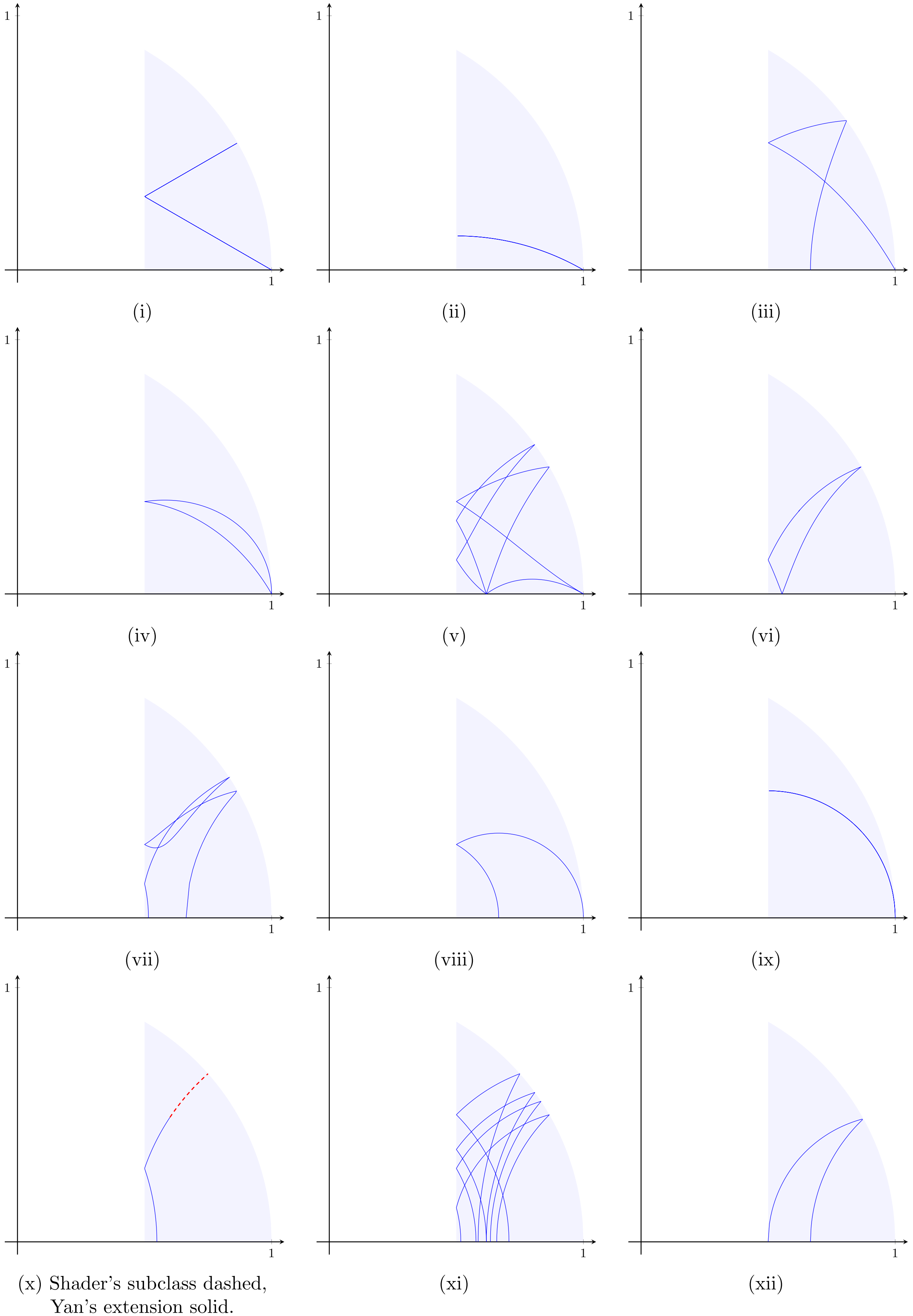

To visualize these classes, we adopt the following useful conventions from~\cite{aich}. Consider a triangle with vertices $A$, $B$, and $C$. After possibly relabeling the vertices, we can assume without loss of generality that $|AB|\ge |AC| \ge |BC|$, where $|\cdot|$ stands for the Euclidean norm. Since all the coloring problems in this paper are invariant under scaling, we can also assume that $|AB|=1$. Finally, by applying a suitable planar motion, we assume that $A=(0,0)$, $B=(1,0)$ and that $C$ lies in the upper half-plane. We call triangles satisfying all these assumptions \textit{normed}. As in Figure \ref{fig:normed-tri}, a normed triangle can now be identified only by a single point $C$ in the bounded region, consisting of all the points $(x,y) \in \R^2$ such that $x \ge \frac{1}{2}$, $y \ge 0$, and $x^2+y^2\le 1$. With these conventions, we visualize classes (i)-(xii) in \Cref{fig:classes}. 

\begin{figure}[htbp]
    \centering
    \resizebox{0.45\linewidth}{!}{
      \begin{tikzpicture}
\begin{axis}[myplot,
    width=8cm,
    height=8cm]
\addplot[name path=f, domain=.5:1, samples=200, red, draw=none] {sqrt(1-x^2)};
\path[name path=axis] (axis cs:0,0) -- (axis cs:1,0) {};
\addplot[
        thick,
        color=blue,
        fill=blue, 
        fill opacity=0.05
    ]
fill between[
  of=f and axis,
  soft clip={domain=.5:1}
];

\coordinate[label={[label distance=3mm]80:$A$}] (A) at (0,0);
\coordinate[label={[label distance=3mm]160:$B$}] (B) at (1,0);
\coordinate[label=above:$C$] (C) at (0.7, 0.5);
\draw[dashed, strongblue, line width=1pt] (A) -- (B) -- (C) --cycle;
\fill (A) circle (2pt);
\fill (B) circle (2pt);
\fill (C) circle (2pt);
\end{axis}
\end{tikzpicture}
    }
    \captionsetup{justification=centering}
    \caption{The shaded area indicates all possible positions of $C$ \\ such that the set $\{A,B,C\}\subset\R^2$ forms a normed triangle.}
    \label{fig:normed-tri}
\end{figure}
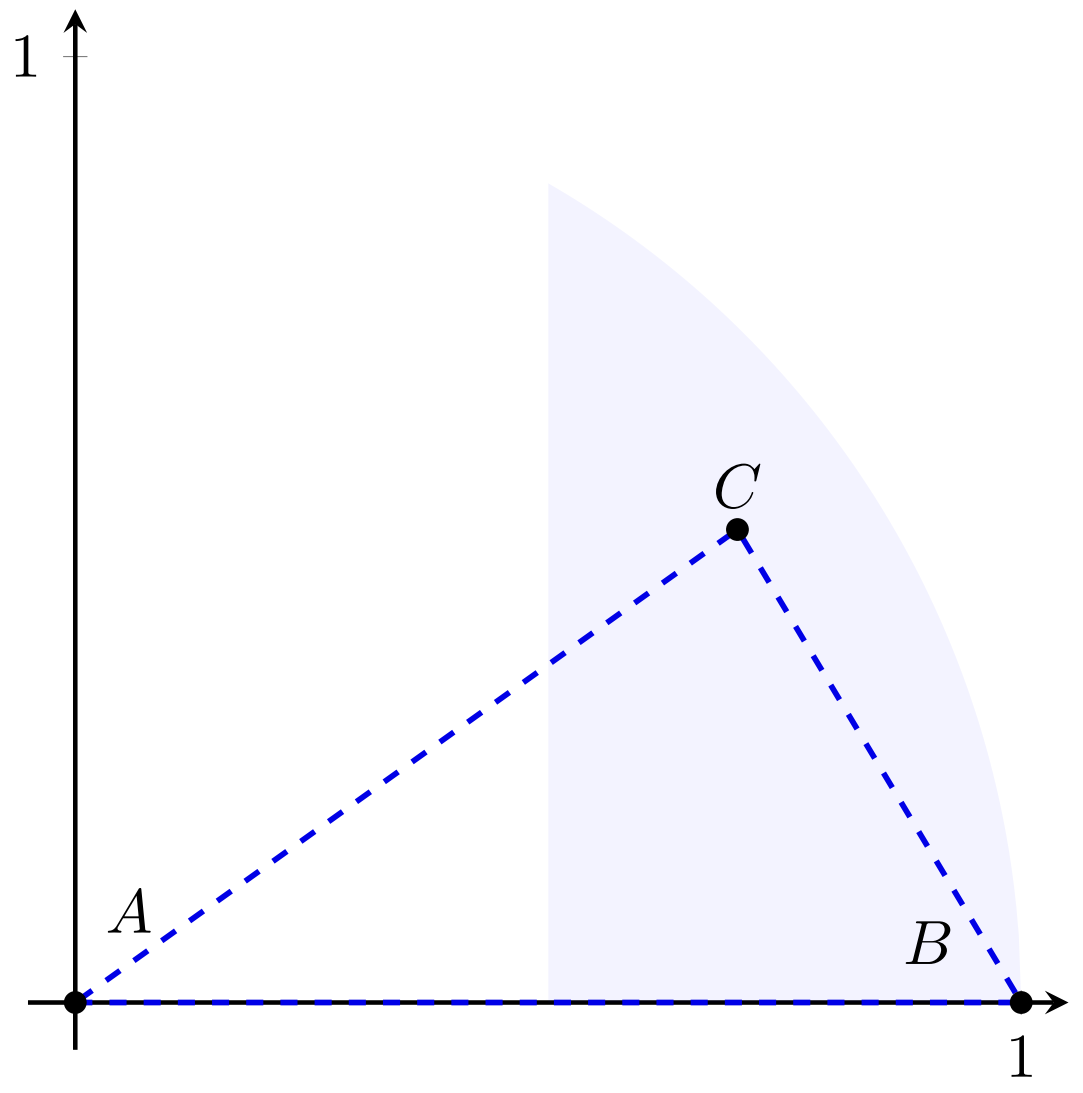

\newpage
\section{Computational approach} \label{sec:comp}

This section contains computer-assisted proofs of \Cref{th:2} and \Cref{th:3}. Both utilize the following standard observation, which we prove for completeness.

\begin{Lemma}
  \label{lem:distI}
  For each $d>0$, every non-monochromatic coloring of the plane contains two points of different colors at distance $d$.
\end{Lemma}

\begin{proof}
    Let $X,Y \in \R^2$ be of different colors. Consider a sequence of points $X=Z_0, Z_1, \ldots, Z_n=Y$ such that $|Z_{i-1}Z_{i}|=d$ for every $1 \le i\le n$. For instance, subdivide the segment $XY$ into subsegments of length at most $2d$, take their endpoints, and between each pair of consecutive endpoints, insert a new point into this sequence at distance $d$ from both. Since the endpoints of this sequence are of different colors, some two consecutive points are of different colors as well.
\end{proof}

All the arguments below take place on finite subsets of the triangular lattice, which we now define. Let $e_1=(1,0)$ and $e_2=(\frac{1}{2}, \frac{\sqrt{3}}{2})$.  For $x,y\in\Z$, let
\( \trip{x}{y} = x\cdot e_1 + y \cdot e_2\), see \Cref{fig:tri_coord}.
The \textit{triangular lattice} $\T$ is the set of all these points:
\[\T \coloneqq \left\{\trip{x}{y}:\, x,y\in\Z\right\}\subset \R^2.\]
For $r\in\N$, we further define 
\[\T_r \coloneqq \left\{\trip{x}{y}\in\T:\, |x|,|y|\le r \right\}\subset \T.\]
The points of $\T_r$ lie in a parallelogram of side lengths $2r$, see \Cref{fig:box}.

\begin{figure}[p]
  \centering
  \begin{minipage}[b]{0.39\textwidth}
    \vspace{0pt}
    \centering
    \resizebox{\linewidth}{!}{%
      \begin{tikzpicture}
        \def\trisize{1.3}
        \trigrid{4}{4}
        \coordinate[label=below:$\trip{0}{0}$] (A) at (\tricoord{0}{1});
        \coordinate (A1) at (\tricoord{1}{1});
        \coordinate (B1) at (\tricoord{0}{2});
        \coordinate[label=below:$\trip{2}{1}$] (C) at (\tricoord{2}{2});
        \draw[vector, strongred] (A) -- (A1) node[midway,above] {$e_1$};
        \draw[vector, strongblue] (A) -- (B1) node[midway,left] {$e_2$};
        \fill (A) circle (2pt);
        \fill (C) circle (2pt);
      \end{tikzpicture}%
    }
    \caption{Triangular lattice $\T$.}
    \label{fig:tri_coord}
  \end{minipage}
  \hfill
  \begin{minipage}[b]{0.45\textwidth}
    \vspace{0pt}
    \centering
    \resizebox{\linewidth}{!}{%
      \begin{tikzpicture}
        \def\trisize{1}
        \trigrid{6}{7}
        \foreach \i in {0,...,4}{
          \foreach \j in {1,...,5}{
            \fill (\tricoord{\i}{\j}) circle (2pt);
          }
        }

        \coordinate[label=below:$\trip{0}{0}$] (A) at (\tricoord{2}{3});
        \fill (\tricoord{3}{3}) circle (2pt);
      \end{tikzpicture}%
    }
    \caption{The bounded subset $\T_2$.}
    \label{fig:box}
  \end{minipage}
\end{figure}

Observe that a positive real occurs as a distance between two points of $\T$ if and only if its square is a \textit{Löschian number}, i.e., of the form $\left|\trip{x}{y}\right|^2  = x^2 + xy + y^2$ for some $x,y \in \Z$. We emphasize that $\T$ not necessarily contain three points forming a triangle $\Tr{a}{b}{c}$ even if all three numbers $a^2$, $b^2$, and $c^2$ are Löschian, as shown by $\Tr{1}{\sqrt{3}}{\sqrt{3}}$. For an explicit characterization of such triangles, see Aebi and Cairns~\cite{aebi}.

\subsection{Proof of Theorem~\ref{th:2}} \label{sec:th:2}

In this approach, we consider the set of triangles $$\mathcal{T} \coloneqq \big\{\{(0,0)_\T,B,C\}: B,C \in \T_4\big\}.$$ Note that there are only 250 distinct triangles in $\mathcal{T}$ up to similarity. For each $\TrO \in \mathcal{T}$, we search for $1 \le d \le 12$ such that every 2-coloring of $\T_6$, in which the points $X=(-6,0)_\T$ and $Y=(-6+d,0)_\T$ receive different colors, contains a monochromatic copy of $\TrO$. To do so, we make use of the z3-solver, a theorem prover from Microsoft Research~\cite{z3}. We claim that if such $d$ exists, then $\R^2 \xrightarrow{2} \TrO$.

Indeed, assume the contrary, namely that there exists a $2$-coloring of the plane containing no monochromatic copies of $\TrO$. In particular, this coloring is non-monochromatic, and thus contains two points of different colors at distance $d$ by \Cref{lem:distI}. By applying a suitable planar motion, we can assume without loss of generality that these points coincide with $X$ and $Y$. Now the 2-coloring induced on $\T_6$ contains no monochromatic copies of $\TrO$ as well while the points $X$ and $Y$ receive different colors in it. However, there is no such 2-coloring by the claim's assumption. This contradiction proves the claim.

The algorithm described above is given in Appendix~\ref{sec:app}. It finds the set $\mathcal{A}_{\text{\upshape M}}$ of 32 new triangles $\TrO$ satisfying $\R^2 \xrightarrow{2} \TrO$, which completes the proof of \Cref{th:2}. It also rediscovers 56 triangles from classes (i)-(xii), see Appendix~\ref{sec:appA}. For the remaining 162 distinct triangles from $\mathcal{T}$, the algorithm does not find the desired $d$, which by itself has no implication on the status of \Cref{con:chromLB} in these special cases.  Since the set $\mathcal{A}_{\text{\upshape M}}$ contains an isosceles triangle $\newT$, \Cref{lem:ratio} now implies that $\R^2 \xrightarrow{2} \Tr{a}{b}{c}$ if $\frac{b}{a}=\sqrt{\frac{7}{3}}$, which proves \Cref{th:1}. \Cref{fig:allMono} visualizes this new class, along with 88 triangles from Appendix~\ref{sec:appA}.

For the sake of completeness, we provide a self-contained proof of \Cref{th:1} in \Cref{sec:4}.

\subsection{Proof of Theorem~\ref{th:3}}

In a very similar fashion to \Cref{sec:th:2}, we again consider the set of triangles $\mathcal{T}$.  For each $\TrO \in \mathcal{T}$, we search for $1 \le d \le 12$ such that every coloring of $\T_6$, in which the points $X=(-6,0)_\T$ and $Y=(-6+d,0)_\T$ receive different colors, regardless of the number of colors used, contains a monochromatic or a rainbow copy of $\TrO$. To do so, we make use of the z3-solver~\cite{z3}. As earlier in \Cref{sec:th:2}, \Cref{lem:distI} implies that if such $d$ exists, then $\R^2 \MRarrows \TrO$.

We again refer to Appendix~\ref{sec:app} for the algorithm described above, as it only differs from the 2-colored case by removing the restriction \lstinline{c < 2} in Line 77. This algorithm finds the set $\mathcal{A}_{\text{\upshape MR}}$ of the first known 36 triangles $\TrO$ satisfying $\R^2 \MRarrows \TrO$, which completes the proof of \Cref{th:3}. Only 21 of these triangles $\TrO$ were known to satisfy the weaker property $\R^2 \xrightarrow{2} \TrO$, see \Cref{fig:foundMix} and Appendix~\ref{sec:appA}.

For a specific right triangle $\TrO=\Tr{1}{\sqrt{3}}{2}$ from $\mathcal{A}_{\text{\upshape MR}}$, we reprove that $\R^2 \MRarrows \TrO$ using a self-contained argument in \Cref{sec:5} to illustrate the additional difficulty arising when more than 2 colors are allowed.

\begin{figure}[htbp]
  \centering
  \begin{minipage}[t]{0.49\textwidth}
    \vspace{0pt}
    \centering
    \resizebox{\linewidth}{!}{
      \begin{tikzpicture}
        \begin{axis}[myplot,
        height=10cm]
          \addplot[name path=f, domain=.5:1, samples=200, red, draw=none] {sqrt(1-x^2)};
          \path[name path=axis] (axis cs:0,0) -- (axis cs:1,0);
          \addplot[
            thick,
            color=blue,
            fill=blue, 
            fill opacity=0.05,
          ]
          fill between[
            of=f and axis,
            soft clip={domain=.5:1}
          ];
          
          \addplot [
            thin,
            blue,
            mark=none,
            draw opacity=\lineopacity,
          ] table {curve_data/30deg.dat};
          \addplot [
            thin,
            blue,
            mark=none,
            draw opacity=\lineopacity,
          ] table {curve_data/150deg.dat};
          \addplot [
            thin,
            blue,
            mark=none,
            draw opacity=\lineopacity,
          ] table {curve_data/a_2a.dat};
          \addplot [
            thin,
            blue,
            mark=none,
            draw opacity=\lineopacity,
          ] table {curve_data/180minusA.dat};
          \addplot [
            thin,
            blue,
            mark=none,
            draw opacity=\lineopacity,
          ] table {curve_data/comp1_1.dat};
          \addplot [
            thin,
            blue,
            mark=none,
            draw opacity=\lineopacity,
          ] table {curve_data/comp1_2.dat};
          \addplot [
            thin,
            blue,
            mark=none,
            draw opacity=\lineopacity,
          ] table {curve_data/comp2.dat};
          \addplot [
            thin,
            blue,
            mark=none,
            draw opacity=\lineopacity,
          ] table {curve_data/isoc.dat};
          \addplot [
            thin,
            blue,
            mark=none,
            draw opacity=\lineopacity,
          ] table {curve_data/ratio_30deg.dat};
          \addplot [
            thin,
            blue,
            mark=none,
            draw opacity=\lineopacity,
          ] table {curve_data/ratio_36deg.dat};
          \addplot [
            thin,
            blue,
            mark=none,
            draw opacity=\lineopacity,
          ] table {curve_data/ratio_90deg.dat};
          \addplot [
            thin,
            blue,
            mark=none,
            draw opacity=\lineopacity,
          ] table {curve_data/ratio_120deg.dat};
          \addplot [
            thin,
            blue,
            mark=none,
            draw opacity=\lineopacity,
          ] table {curve_data/90deg.dat};
          \addplot [
            thin,
            blue,
            mark=none,
            draw opacity=\lineopacity,
          ] table {curve_data/aSquaredPlus2BSquared.dat};
          \addplot [
            thin,
            blue,
            mark=none,
            draw opacity=\lineopacity,
          ] table {curve_data/bach.dat};
          \addplot [
            thin,
            blue,
            mark=none,
            draw opacity=\lineopacity,
          ] table {curve_data/degen.dat};
          \addplot [
            thin,
            dash pattern=on 2pt off 0.5pt,
            red,
            mark=none,
            draw opacity=0.7,
          ] table {curve_data/new.dat};
          \addplot[
            only marks,
            strongblue,
            mark size=2pt,
            mark=+,
          ] table {point_data/knownMono.dat};
          \addplot[
            only marks,
            strongred,
            mark=x,
            mark size=2pt,
          ] table {point_data/unknownMono.dat};
          
        \end{axis}
      \end{tikzpicture}
    }
    \caption{The points identifying all normed triangles $\TrO$ known to satisfy  $\R^2\rArrows{2}\TrO$. Blue lines correspond to classes (i)-(xii), while the red dashed line corresponds to the new class from \Cref{th:1}. Individual points correspond to the triangles from  Appendix~\ref{sec:appA}, the blue ones (+) lie on the blue, previously known lines, while the red ones ($\times$) correspond to the new results from \Cref{th:2}.} 
    \label{fig:allMono}
  \end{minipage}
  \hfill
  \begin{minipage}[t]{0.49\textwidth}
    \vspace{0pt}
    \centering
    \resizebox{\linewidth}{!}{
    \begin{tikzpicture}
        \begin{axis}[myplot,
        height=10cm]
          \addplot[name path=f, domain=.5:1, samples=200, red, draw=none] {sqrt(1-x^2)};
          \path[name path=axis] (axis cs:0,0) -- (axis cs:1,0);
          \addplot[
            thick,
            color=blue,
            fill=blue,
            fill opacity=0.05,
          ]
          fill between[
            of=f and axis,
            soft clip={domain=.5:1}
          ];
          \addplot [
            thin,
            dash pattern=on 1.0pt off 2.0pt,
            blue,
            mark=none,
            draw opacity=\lineopacity,
          ] table {curve_data/90deg.dat};
          \addplot[
            only marks,
            strongred,
            mark=x,
            mark size=2pt,
          ] table {point_data/foundMix_unknownMono.dat};
          \addplot[
            only marks,
            strongblue,
            mark=+,
            mark size=2pt,
          ] table {point_data/foundMix_knownMono.dat};
        \end{axis}
      \end{tikzpicture}
    }
    \caption{The points identifying all 36 normed triangles $\TrO$ known to satisfy  $\R^2 \MRarrows \TrO$, see \Cref{th:3}. The blue ones (+) lie on the blue lines from \Cref{fig:allMono}, while the red ones ($\times$) correspond to the triangles $\TrO$ for which even the weaker property $\R^2 \xrightarrow{2} \TrO$ is new.
    The dashed blue line represents right triangles; the points corresponding to acute and obtuse triangles lie above and below it, respectively.} 
    \label{fig:foundMix}
  \end{minipage}
\end{figure}

\newpage
\enlargethispage{3mm}
\section{The special case \texorpdfstring {$\TrO=\newT$} {T = T(√3, √7, √7)} of Theorem~\ref{th:2}} \label{sec:4}

In this section, we reprove the following result without computer assistance. 

\begin{Proposition}\label{thm:newT}
  We have $\R^2\rArrows{2}\newT$.
\end{Proposition}

Since $\newT$ is isosceles, \Cref{lem:ratio} and \Cref{thm:newT} together yield \Cref{th:1}.

\begin{proof} [Proof of \Cref{th:1}, assuming \Cref{thm:newT}]
    Let $\TrO=\Tr{a}{b}{c}$ be a triangle such that $\frac{b}{a}=\sqrt{\frac{7}{3}}$. Let $\TrO'$ be a copy of $\TrO$ scaled by the factor $\frac{\sqrt{3}}{a}$, namely $\TrO'=\Tr{\sqrt{3}}{\sqrt{7}}{c\frac{\sqrt{3}}{a}}$. \Cref{lem:ratio} and \Cref{thm:newT} imply that $\R^2 \xrightarrow{2} \TrO'$. Now a standard scaling argument shows that $\R^2 \xrightarrow{2} \TrO$, as desired.
    Indeed, assume the contrary, namely that there exists a $2$-coloring of the plane containing no monochromatic copies of $\TrO$. Then the auxiliary 2-coloring scaled by the factor $\frac{\sqrt{3}}{a}$ contains no monochromatic copies of the scaled triangle $\TrO'$, a contradiction.
\end{proof}

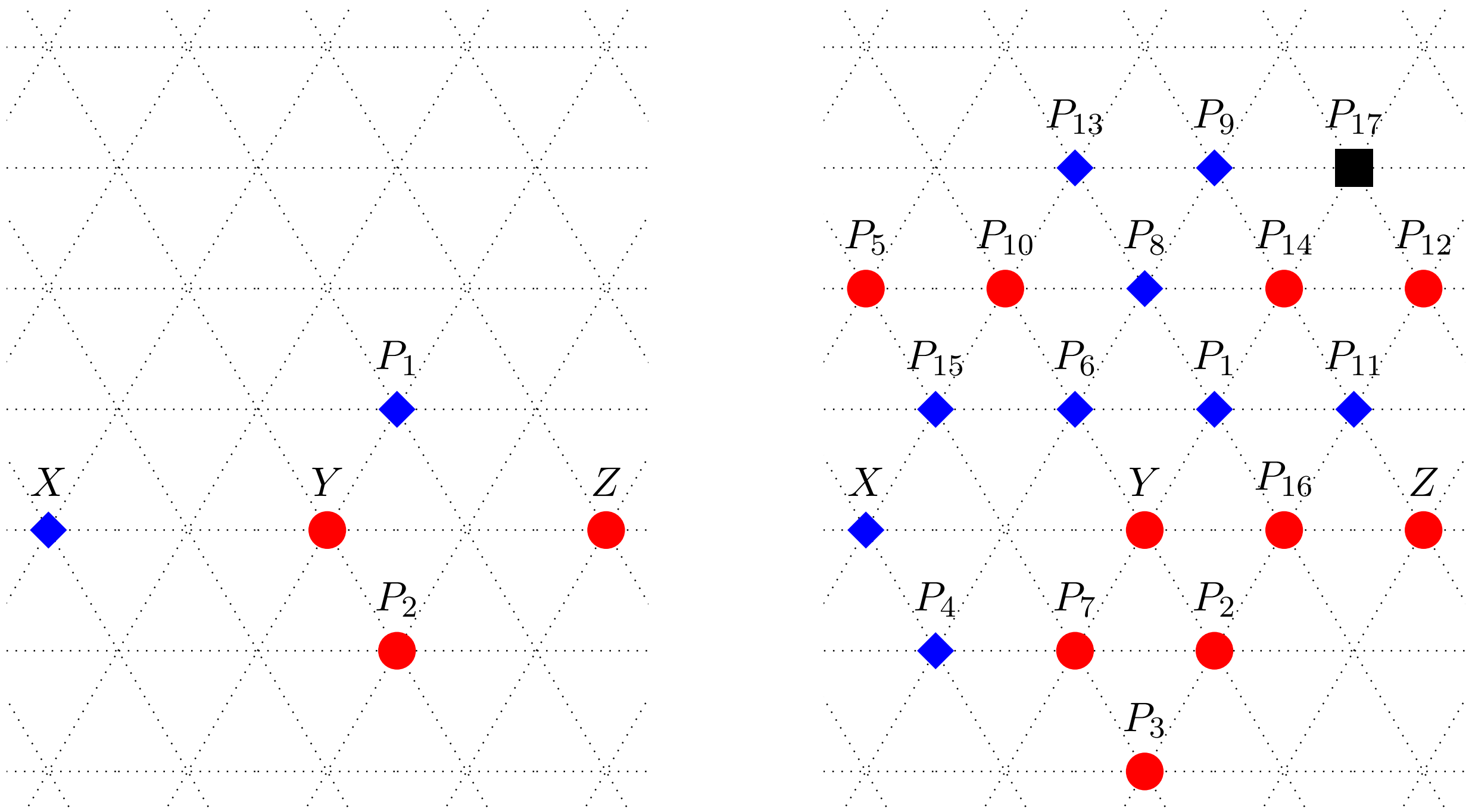
\begin{figure}[htbp]
  \begin{subfigure}[t]{0.45\textwidth}
    \centering
    \resizebox{\linewidth}{!}{
      \begin{tikzpicture}
        \def\trisize{1.3}
        \trigrid{6}{4}
        
        \node[label=$X$, diamond, fill=blue, minimum size=10pt, inner sep=0pt] at (\tricoord{-1}{2}) {};
        \node[label=$Y$, circle, fill=red, minimum size=10pt, inner sep=0pt] at (\tricoord{1}{2}) {};
        \node[label=$Z$, circle, fill=red, minimum size=10pt, inner sep=0pt] at (\tricoord{3}{2}) {};

        \node[label=$P_2$, circle, fill=red, minimum size=10pt, inner sep=0pt] at (\tricoord{2}{1}) {};
        \node[label=$P_1$, diamond, fill=blue, minimum size=10pt, inner sep=0pt] at (\tricoord{1}{3}) {};
      \end{tikzpicture}
    }
    \caption{Initial construction.}
    \label{fig:newT1}
  \end{subfigure}
  \hfill
  \begin{subfigure}[t]{0.45\textwidth}
    \centering
    \resizebox{\linewidth}{!}{
      \begin{tikzpicture}
        \def\trisize{1.3}
        \trigrid{6}{4}

        \node[label=$X$, diamond, fill=blue, minimum size=10pt, inner sep=0pt] at (\tricoord{-1}{2}) {};
        
        \node[label=$Y$, circle, fill=red, minimum size=10pt, inner sep=0pt] at (\tricoord{1}{2}) {};
        \node[label=$Z$, circle, fill=red, minimum size=10pt, inner sep=0pt] at (\tricoord{3}{2}) {};

        \node[label=$P_1$, diamond, fill=blue, minimum size=10pt, inner sep=0pt] at (\tricoord{1}{3}) {};
        \node[label=$P_2$, circle, fill=red, minimum size=10pt, inner sep=0pt] at (\tricoord{2}{1}) {};
        \node[label=$P_3$, circle, fill=red, minimum size=10pt, inner sep=0pt] at (\tricoord{2}{0}) {};
        \node[label=$P_4$, diamond, fill=blue, minimum size=10pt, inner sep=0pt] at (\tricoord{0}{1}) {};
        \node[label=$P_5$, circle, fill=red, minimum size=10pt, inner sep=0pt] at (\tricoord{-2}{4}) {};
        \node[label=$P_6$, diamond, fill=blue, minimum size=10pt, inner sep=0pt] at (\tricoord{0}{3}) {};
        \node[label=$P_7$, circle, fill=red, minimum size=10pt, inner sep=0pt] at (\tricoord{1}{1}) {};
        \node[label=$P_8$, diamond, fill=blue, minimum size=10pt, inner sep=0pt] at (\tricoord{0}{4}) {};
        \node[label=$P_9$, diamond, fill=blue, minimum size=10pt, inner sep=0pt] at (\tricoord{0}{5}) {};
        \node[label=$P_{10}$, circle, fill=red, minimum size=10pt, inner sep=0pt] at (\tricoord{-1}{4}) {};
        \node[label=$P_{11}$, diamond, fill=blue, minimum size=10pt, inner sep=0pt] at (\tricoord{2}{3}) {};
        \node[label=$P_{12}$, circle, fill=red, minimum size=10pt, inner sep=0pt] at (\tricoord{2}{4}) {};
        \node[label=$P_{13}$, diamond, fill=blue, minimum size=10pt, inner sep=0pt] at (\tricoord{-1}{5}) {};
        \node[label=$P_{14}$, circle, fill=red, minimum size=10pt, inner sep=0pt] at (\tricoord{1}{4}) {};
        \node[label=$P_{15}$, diamond, fill=blue, minimum size=10pt, inner sep=0pt] at (\tricoord{-1}{3}) {};
        \node[label=$P_{16}$, circle, fill=red, minimum size=10pt, inner sep=0pt] at (\tricoord{2}{2}) {};

        \node[label=$P_{17}$, rectangle, fill=black, minimum size=10pt, inner sep=0pt] at (\tricoord{1}{5}) {};
      \end{tikzpicture}
    }
    \caption{Auxiliary points.}
    \label{fig:newTall}
  \end{subfigure}
  \par
  \caption{Constructions in the proof of \Cref{thm:newT}.}
  \label{fig:newTproof}
\end{figure}

\begin{proof}[Proof of \Cref{thm:newT}]
  Assume the contrary, namely that there exists a $2$-coloring of the plane in red and blue containing no monochromatic copies of $\newT$. Recall that this coloring also contains no monochromatic copies of $\TrI{\sqrt{3}}$ and $\TrI{\sqrt{7}}$ by \Cref{theorem:equi}. In particular, this coloring is non-monochromatic, and thus contains two points of different colors at distance $4$ by \Cref{lem:distI}. By applying a suitable planar motion, we can assume without loss of generality that these points are blue $X=\trip{0}{0}$ and red $Z=\trip{4}{0}$. By symmetry, we can also assume that their midpoint $Y=\trip{2}{0}$ is red.
  The points $P_1=\trip{2}{1}$ and $P_2=\trip{3}{-1}$ must be of different colors, as otherwise either the copy $\{X,P_1,P_2\}$ of $\newT$ would be blue or the copy $\{Z,P_1,P_2\}$ of $\TrI{\sqrt{3}}$ would be red. By symmetry again, assume that $P_1$ is blue and $P_2$ is red. This initial construction is illustrated in \Cref{fig:newT1}.
  
  Now we reconstruct the colors of some other points of the triangular lattice one by one until we find the desired contradiction. To shorten the process of describing the construction, we use the following table. It names further auxiliary points and justifies their reconstructed colors.

\vspace{3mm}
\begin{tabular}{llll}
Name & Coordinates & Color & Justification (Avoiding monochromatic triangle) \\ \hline
 $P_3$     & $\trip{3}{-2}$ & Red   & $\TrI{\sqrt{7}} : \{X,P_1,P_3\}$ \\
 $P_4$     & $\trip{1}{-1}$ & Blue  & $\TrI{\sqrt{3}} : \{Y,P_3,P_4\}$ \\
 $P_5$     & $\trip{-1}{2}$ & Red   & $\TrI{\sqrt{7}} : \{P_1,P_4,P_5\}$ \\
 $P_6$     & $\trip{1}{1}$  & Blue  & $\TrI{\sqrt{7}} : \{Z,P_3,P_6\}$ \\
 $P_7$     & $\trip{2}{-1}$ & Red   & $\TrI{\sqrt{3}} : \{X,P_6,P_7\}$ \\
 $P_8$     & $\trip{1}{2}$  & Blue  & $\TrI{\sqrt{7}} : \{Z,P_7,P_8\}$ \\
 $P_9$     & $\trip{1}{3}$  & Blue  & $\TrI{\sqrt{7}} : \{Y,P_5,P_9\}$ \\
 $P_{10}$   & $\trip{0}{2}$  & Red   & $\TrI{\sqrt{3}} : \{P_1,P_9,P_{10}\}$ \\ 
 $P_{11}$   & $\trip{3}{1}$  & Blue  & $\TrI{\sqrt{7}} : \{P_7,P_{10},P_{11}\}$ \\
 $P_{12}$   & $\trip{3}{2}$  & Red   & $\TrI{\sqrt{3}} : \{P_1,P_9,P_{12}\}$ \\
 $P_{13}$   & $\trip{0}{3}$  & Blue  & $\TrI{\sqrt{7}} : \{Y,P_{12},P_{13}\}$ \\ 
 $P_{14}$   & $\trip{2}{2}$  & Red   & $\TrI{\sqrt{3}} : \{P_6,P_{13},P_{14}\}$ \\ 
 $P_{15}$   & $\trip{0}{1}$  & Blue  & $\TrI{\sqrt{7}} : \{P_2,P_{14},P_{15}\}$ \\ 
 $P_{16}$   & $\trip{3}{0}$  & Red   & $\TrI{\sqrt{7}} : \{P_9,P_{15},P_{16}\}$ \\
\hline
\end{tabular}
\vspace{4mm}

The resulting configuration is illustrated in Figure \ref{fig:newTall}. Now, regardless of the color of $P_{17}=\trip{2}{3}$, either the copy $\{P_{16},P_{10},P_{17}\}$ of $\TrI{\sqrt{7}}$ is red, or the copy $\{P_{11},P_8,P_{17}\}$ of $\TrI{\sqrt{3}}$ is blue. This contradiction shows that $\R^2\rArrows{2}\newT$, as desired.
\end{proof}

\section{The special case \texorpdfstring {$\TrO=\Tr{1}{\sqrt{3}}{2}$} {T = T(1, √3, 2)} of Theorem~\ref{th:3}} \label{sec:5}

In this section, we reprove the following result without computer assistance. 

\begin{Proposition}\label{thm:306090}
  We have $\R^2 \MRarrows \Tr{1}{\sqrt{3}}{2}$.
\end{Proposition}

The goal of this is twofold. First, to confirm \Cref{con:MR} for at least one specific triangle with a self-contained argument. Second, to illustrate the additional difficulty arising when more than 2 colors are allowed. Indeed, showing that $\R^2 \xrightarrow{2} \TrO$ is perhaps easier for $\TrO=\Tr{1}{\sqrt{3}}{2}$ than for any other triangle. The next argument showing that $\R^2 \xrightarrow{2} \Tr{1}{\sqrt{3}}{2}$ is due to Soifer~\cite{soifer} and it takes only one paragraph.

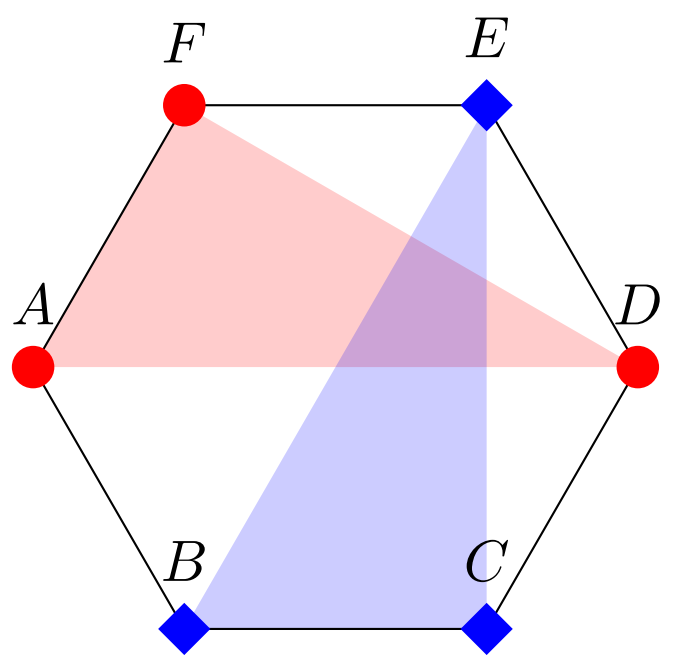
\begin{figure}[htbp]
  \centering
  \begin{tikzpicture}
    \def\trisize{2}

\coordinate (A) at (\tricoord{-1}{0});
\coordinate (B) at (\tricoord{0}{-1});
\coordinate (C) at (\tricoord{1}{-1});
\coordinate (D) at (\tricoord{1}{0});
\coordinate (E) at (\tricoord{0}{1});
\coordinate (F) at (\tricoord{-1}{1});

\fill[red, fill opacity=0.2] (A) -- (D) -- (F) -- cycle;
\fill[blue, fill opacity=0.2] (E) -- (B) -- (C) -- cycle;

\draw (A)--(B)--(C)--(D)--(E)--(F)--cycle;

\node[label=$A$, circle, fill=red, minimum size=8pt, inner sep=0pt] at (A) {};
\node[label=$D$, circle, fill=red, minimum size=8pt, inner sep=0pt] at (D) {};
\node[label=$B$, diamond, fill=blue, minimum size=10pt, inner sep=0pt] at (B) {};
\node[label=$C$, diamond, fill=blue, minimum size=10pt, inner sep=0pt] at (C) {};
\node[label=$E$, diamond, fill=blue, minimum size=10pt, inner sep=0pt] at (E) {};
\node[label=$F$, circle, fill=red, minimum size=8pt, inner sep=0pt] at (F) {};
\end{tikzpicture}
\caption{Construction used to prove that $\R^2 \xrightarrow{2} \Tr{1}{\sqrt{3}}{2}$.}
\label{fig:mono306090}
\end{figure}

Consider an arbitrary $2$-coloring of the plane in red and blue. By the pigeonhole principle, among three vertices at pairwise distances 2, some two receive the same color, say red. Label them $A$ and $D$ and consider the regular hexagon $ABCDEF$ with diameter $AD$, see \Cref{fig:mono306090}. If at least one of the vertices $B,C,E,F$ is red, say $F$, then $\{A,D,F\}$ is a red copy of $\Tr{1}{\sqrt{3}}{2}$. Otherwise, if all four of them are blue, then, say, $\{B,C,E\}$ is a blue copy of $\Tr{1}{\sqrt{3}}{2}$. In all cases, we find a monochromatic copy of $\Tr{1}{\sqrt{3}}{2}$ and thus $\R^2 \xrightarrow{2} \Tr{1}{\sqrt{3}}{2}$.

\newpage
In sharp contrast with Soifer's elegant argument above, our proof of \Cref{thm:306090} below is much more involved and relies on the following result.

\begin{Lemma}
  \label{lemma:gadget}
  If a coloring of the plane contains two points of the same color at distance 1,
  then it also contains a monochromatic or a rainbow copy of $\Tr{1}{\sqrt{3}}{2}$. 
\end{Lemma}

\begin{figure}[b]
  \begin{subfigure}[b]{0.45\textwidth}
    \centering
    \resizebox{0.5\linewidth}{!}{
      \begin{tikzpicture}
        \def\trisize{1.3}
        \trigrid{3}{2}
        \node[label=$X$, diamond, fill=blue, minimum size=10pt, inner sep=0pt] at (\tricoord{0}{0}) {};
        \node[label=$Y$, diamond, fill=blue, minimum size=10pt, inner sep=0pt] at (\tricoord{1}{0}) {};
        \node[label=$Z$, circle, fill=red, minimum size=10pt, inner sep=0pt] at (\tricoord{0}{2}) {};
        \node[label=$P_1$, circle, fill=red, minimum size=10pt, inner sep=0pt] at (\tricoord{-1}{2}) {};
        \node[label=$P_2$, circle, fill=red, minimum size=10pt, inner sep=0pt] at (\tricoord{1}{2}) {};
        \node[label=$P_3$, diamond, fill=blue, minimum size=10pt, inner sep=0pt] at (\tricoord{2}{0}) {};
        \node[label=$P_4$, rectangle, fill=black, minimum size=10pt, inner sep=0pt] at (\tricoord{0}{1}) {};
      \end{tikzpicture}}
    \caption{Case 1.}
    \label{fig:gadgetB}
  \end{subfigure}
  \hfill
  \begin{subfigure}[b]{0.45\textwidth}
    \centering
    \resizebox{0.5\linewidth}{!}{
      \begin{tikzpicture}
        \def\trisize{1.3}
        \trigrid{3}{2}
        \node[label=$X$, diamond, fill=blue, minimum size=10pt, inner sep=0pt] at (\tricoord{0}{0}) {};
        \node[label=$Y$, diamond, fill=blue, minimum size=10pt, inner sep=0pt] at (\tricoord{1}{0}) {};
        \node[label=$Z$, circle, fill=red, minimum size=10pt, inner sep=0pt] at (\tricoord{0}{2}) {};
        \node[label=$P_1$, circle, fill=red, minimum size=10pt, inner sep=0pt] at (\tricoord{-1}{2}) {};
        \node[label=$P_2$, diamond, fill=blue, minimum size=10pt, inner sep=0pt] at (\tricoord{1}{2}) {};
        \node[label=$P_3$, circle, fill=red, minimum size=10pt, inner sep=0pt] at (\tricoord{2}{0}) {};
        \node[label=$P_4$, rectangle, fill=black, minimum size=10pt, inner sep=0pt] at (\tricoord{0}{1}) {};
      \end{tikzpicture}}
    \caption{Case 2.}
    \label{fig:gadgetA}
  \end{subfigure}
  \par
  \caption{Constructions used in the proof of \Cref{lemma:gadget}.}
  \label{fig:gadget}
\end{figure}
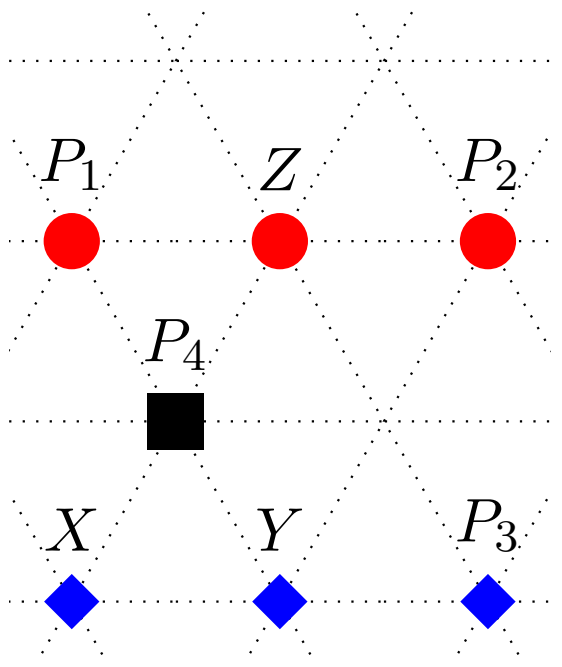
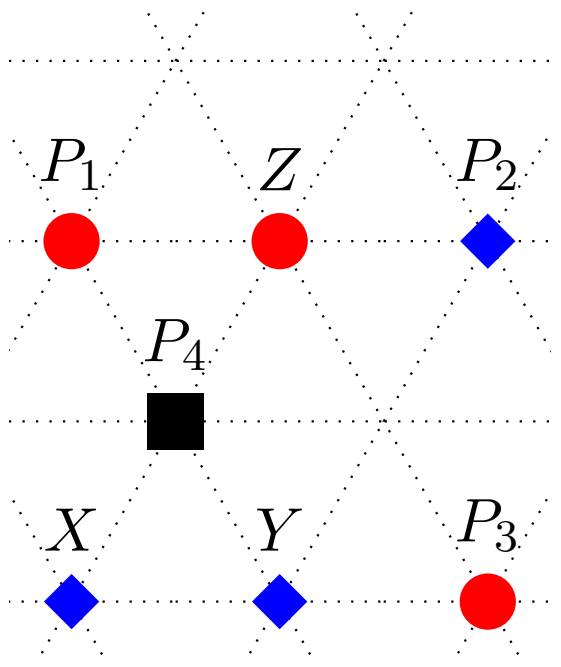

\begin{proof}
  Without loss of generality, let the points of the same color at distance 1 be blue $X=\trip{0}{0}$ and $Y=\trip{1}{0}$. If $Z=\trip{0}{2}$ is blue, then $\{X,Y,Z\}$ is the desired  blue copy of $\Tr{1}{\sqrt{3}}{2}$. So let $Z$ be of different color, say red.

  Now we reconstruct the colors of some other lattice points while avoiding monochromatic and rainbow copies of $\Tr{1}{\sqrt{3}}{2}$ until we get the desired conclusion that this is impossible.
  First, the point $P_1=\trip{-1}{2}$ must be red to avoid the blue copy $\{X,Y,P_1\}$ and the rainbow copy $\{Y,Z,P_1\}$. Next, the point $P_2=\trip{1}{2}$ must be either red or blue to avoid the rainbow copy $\{Y,Z,P_2\}$. We consider these two cases separately.

\begin{description}[
  font=\normalfont\itshape
]
  \item[Case 1: $P_2$ is red.]
    The point $P_3=\trip{2}{0}$ must be blue to avoid the red copy $\{Z,P_2,P_3\}$ and the rainbow copy $\{Y,P_2,P_3\}$. Now, regardless of the color of $P_{4}=\trip{0}{1}$, the copy $\{X,P_3,P_4\}$ is blue, or the copy $\{P_1,P_2,P_4\}$ is red, or the copy $\{Y,P_2,P_4\}$ is rainbow, see \Cref{fig:gadgetB}. 
  \item[Case 2: $P_2$ is blue.]
    The point $P_3=\trip{2}{0}$ must be red to avoid the blue copy $\{Y,P_2,P_3\}$ and the rainbow copy $\{Z,P_2,P_3\}$. Now, regardless of the color of $P_{4}=\trip{0}{1}$, the copy $\{Y,P_2,P_4\}$ is blue, or the copy $\{Z,P_3,P_4\}$ is red, or the copy $\{P_1,P_2,P_4\}$ is rainbow, see \Cref{fig:gadgetA}. \qedhere
  
\end{description}
\end{proof}
\enlargethispage{4mm}

\vspace{-4mm}

\begin{figure}[htbp]
    \centering
  \begin{subfigure}[b]{0.45\textwidth}
    \centering
    \resizebox{\linewidth}{!}{
      \begin{tikzpicture}
    \def\trisize{1.2}
    \trigrid{5}{6}
    \draw (\tricoord{1}{2}) -- (\tricoord{1}{3}) -- (\tricoord{3}{2}) -- cycle;
    \node[label=below:$X$] at (\tricoord{1}{2}) {};
    \node[label=below:$Y$] at (\tricoord{3}{2}) {};
    \node[label=below:$Z$] at (\tricoord{1}{3}) {};
    \node[label=below:$P$] at (\tricoord{2}{2}) {};
    \node[label=below:$Q$] at (\tricoord{2}{3}) {};
    \node[label=below:$R$] at (\tricoord{1}{4}) {};
    \node[label=below:$S$] at (\tricoord{2}{4}) {};
    \def\xoff{0}
    \def\yoff{0}
    \foreach \i in {-1,...,5}{
    \foreach \j in {0,...,5}{
      \def\ci{2 * \i}
      \def\cj{2 * \j}
      \node[regular polygon, fill=yellow, minimum size=6pt, inner sep=0pt] at (\tricoord{\ci+\xoff}{\cj+\yoff}) {};
    }
    }
    \def\xoff{1}
    \def\yoff{0}
    \foreach \i in {-1,...,4}{
    \foreach \j in {0,...,5}{
      \def\ci{2 * \i}
      \def\cj{2 * \j}
      \node[circle, fill=red, minimum size=6pt, inner sep=0pt] at (\tricoord{\ci+\xoff}{\cj+\yoff}) {};
    }
    }
    \def\xoff{1}
    \def\yoff{1}
    \foreach \i in {-1,...,4}{
    \foreach \j in {0,...,5}{
      \def\ci{2 * \i}
      \def\cj{2 * \j}
      \node[diamond, fill=blue, minimum size=8pt, inner sep=0pt] at (\tricoord{\ci+\xoff}{\cj+\yoff}) {};
    }
    }
    \def\xoff{0}
    \def\yoff{1}
    \foreach \i in {-1,...,4}{
    \foreach \j in {0,...,5}{
      \def\ci{2 * \i}
      \def\cj{2 * \j}
      \node[rectangle, fill=green, minimum size=6pt, inner sep=0pt] at (\tricoord{\ci+\xoff}{\cj+\yoff}) {};
    }
    }
  \end{tikzpicture}}
    \caption{Case 1: 2-periodic 4-coloring.}
    \label{fig:per2}
  \end{subfigure}
  \hfill
  \begin{subfigure}[b]{0.45\textwidth}
    \centering
    \resizebox{\linewidth}{!}{
      \begin{tikzpicture}
    \def\trisize{1.2}
    \trigrid{5}{6}
    \node[label=below:$X$] at (\tricoord{1}{2}) {};
    \node[label=below:$Y$] at (\tricoord{3}{2}) {};
    \node[label=below:$Z$] at (\tricoord{1}{3}) {};
    \node[label=below:$P$] at (\tricoord{2}{2}) {};
    \node[label=below:$Q$] at (\tricoord{2}{3}) {};
    \node[label=below:$R$] at (\tricoord{1}{4}) {};
    \node[label=below:$S$] at (\tricoord{2}{4}) {};
    \draw (\tricoord{1}{2}) -- (\tricoord{1}{3}) -- (\tricoord{3}{2}) -- cycle;
    \def\xoff{0}
    \def\yoff{0}
    \foreach \i in {-1,...,5}{
    \foreach \j in {-2,...,5}{
      \def\ci{\i - \j}
      \def\cj{\i + 2*\j}
      \node[rectangle, fill=green, minimum size=6pt, inner sep=0pt] at (\tricoord{\ci+\xoff}{\cj+\yoff}) {};
    }
    }
    \def\xoff{1}
    \def\yoff{0}
    \foreach \i in {-1,...,5}{
    \foreach \j in {-2,...,5}{
      \def\ci{\i - \j}
      \def\cj{\i + 2*\j}
      \node[circle, fill=red, minimum size=6pt, inner sep=0pt] at (\tricoord{\ci+\xoff}{\cj+\yoff}) {};
    }
    }
    \def\xoff{2}
    \def\yoff{0}
    \foreach \i in {-1,...,5}{
    \foreach \j in {-2,...,5}{
      \def\ci{\i - \j}
      \def\cj{\i + 2*\j}
      \node[diamond, fill=blue, minimum size=8pt, inner sep=0pt] at (\tricoord{\ci+\xoff}{\cj+\yoff}) {};
    }
    }
  \end{tikzpicture}}
    \caption{Case 2: 3-periodic 3-coloring.}
    \label{fig:per3}
  \end{subfigure}
  \par
  \caption{The only two colorings of $\T$ avoiding monochromatic and rainbow copies of $\Tr{1}{\sqrt{3}}{2}$.}
  \label{fig:periodic}
\end{figure}
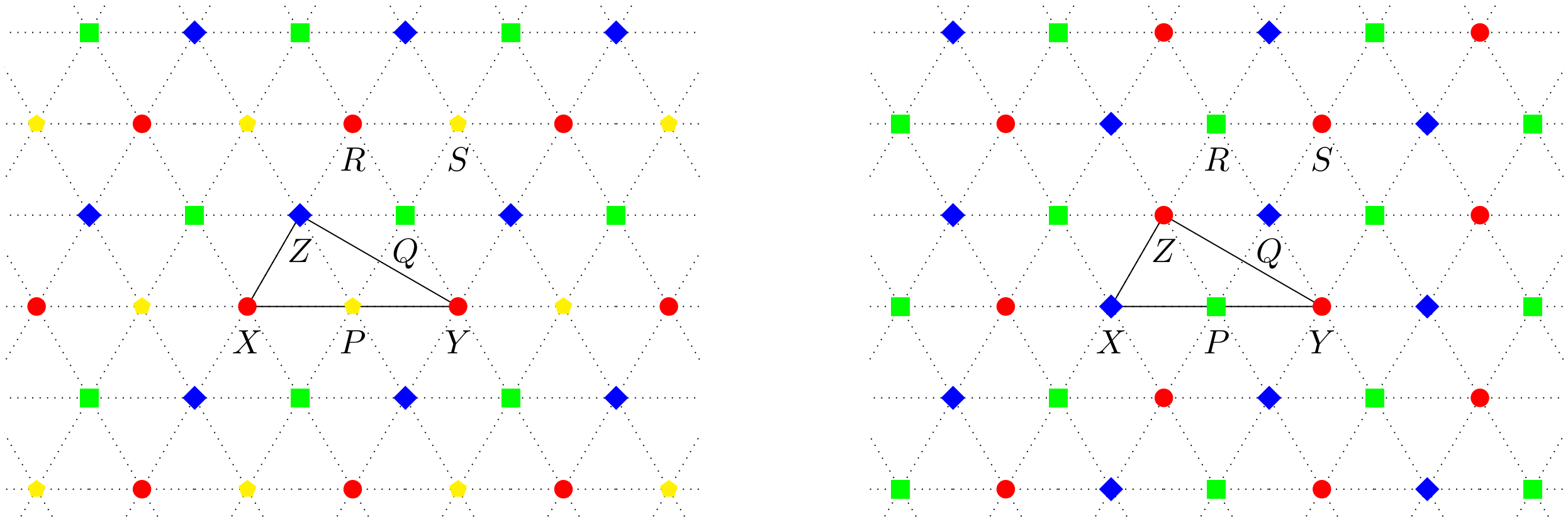

\begin{proof}[Proof of \Cref{thm:306090}]
  Assume the contrary, namely that there exists a coloring of the plane containing no monochromatic and no rainbow copies of $\Tr{1}{\sqrt{3}}{2}$. In particular, the points $X=\trip{0}{0}$, $Y=\trip{2}{0}$ and $Z=\trip{0}{1}$ receive exactly two colors, say two points are red and the third one is blue. Moreover, every two points at distance 1 receive different colors by \Cref{lemma:gadget}. Thus $X$ and $Z$ cannot be both red, so one of them is blue. We consider these two cases separately.
  
  \begin{description}[
  font=\normalfont\itshape]
    \item[Case 1: $Z$ is blue.]
      The points $P=\trip{1}{0}$ and $Q=\trip{1}{1}$ are at distance 1 from each other and from both red $Y$ and blue $Z$. By \Cref{lemma:gadget}, they receive two new colors, say, yellow and green, respectively. The 4-colored parallelogram $XPQZ$ already uniquely determines the colors of the two adjacent points $R=\trip{0}{2}$ and $S=\trip{1}{2}$, see \Cref{fig:per2}. Indeed, the point $R$ must be either red or yellow to avoid the rainbow copy $\{X,P,R\}$ and also cannot be yellow to avoid the rainbow copy $\{X,Q,R\}$, so $R$ is red. Similarly, the point $S$ must be yellow to avoid the rainbow copies $\{Z,P,S\}$ and $\{P,R,S\}$. Thus the parallelogram $ZQSR$, obtained by translating $XPQZ$ by the vector $\trip{0}{1}$, is also 4-colored with the roles of the colors red and yellow switched with those of blue and green, respectively. By reflecting all the points in this argument about the line $XQ$, we reach a similar conclusion about the coloring of the parallelogram, obtained by translating $XPQZ$ by the vector $\trip{1}{0}$. Moreover, this argument applies verbatim not only to $XPQZ$ but also to all its 4-colored translates. Since $\trip{1}{0}$ and $\trip{0}{1}$ form a basis of $\T$, the coloring of $\T$ is now uniquely determined, see \Cref{fig:per2}. In particular, the vector $\trip{0}{2}$ is one of its periods, and thus the point $\trip{0}{6}$ is red. 
    \item[Case 2: $X$ is blue.]
      The point $P=\trip{1}{0}$ is at distance 1 from both blue $X$ and red $Y$, so by \Cref{lemma:gadget}, it receives a new colors, say, green. Similarly, the point $Q=\trip{1}{1}$ cannot be red as it is at distance 1 from red $Y$. Thus $Q$ must be blue to avoid the rainbow copy $\{X,Y,Q\}$ of $\Tr{1}{\sqrt{3}}{2}$. As in Case~1, the 3-colored parallelogram $XPQZ$ with a monochromatic diameter $XQ$ uniquely determines the colors of the two adjacent points $R=\trip{0}{2}$ and $S=\trip{1}{2}$, see \Cref{fig:per3}. Indeed, the point $R$ must be green to avoid the blue copy $\{X,Q,R\}$ and the rainbow copy $\{X,P,R\}$. Similarly, the point $S$ must be red avoid the green copy $\{P,R,S\}$ and the rainbow copy $\{P,Z,S\}$. Thus the parallelogram $ZQSR$, obtained by translating $XPQZ$ by the vector $\trip{0}{1}$, is also 3-colored, with a monochromatic diameter $ZS$ and the roles of the colors cyclically permuted. By reflecting all the points in this argument about the line $XQ$, we reach a similar conclusion about the coloring of the parallelogram, obtained by translating $XPQZ$ by the vector $\trip{1}{0}$. Moreover, this argument applies verbatim not only to $XPQZ$ but also to all its 3-colored translates with monochromatic diameters. As in Case~1, the coloring of $\T$ is now uniquely determined, see \Cref{fig:per3}. In particular, the vector $\trip{0}{3}$ is one of its periods, and thus the point $\trip{0}{6}$ is blue. 
  \end{description}
  
  We have shown that, in both cases, the points $X$ and $X+\trip{0}{6}$ receive the same color. By applying an arbitrary planar motion, we conclude that every two points at distance 6 must receive the same color. Now \Cref{lem:distI} implies that the coloring is monochromatic and, in particular, contains a monochromatic copy of $\Tr{1}{\sqrt{3}}{2}$, a contradiction.
\end{proof}

\newpage
\section{Concluding remarks} \label{sec:concl}

In this paper, we confirmed \Cref{con:chromLB} for a new class of triangles by considering finite configurations on the triangular lattice. One can certainly find many new such classes by considering larger subsets or different lattices, e.g., the \textit{Moser lattice} used recently to solve several related problems, see~\cite{ambrus2024density, currier, ducz2026note, ducz2026unit, engel2025diverse, matolcsi2023fractional, moore}. However, it seems that different ideas are needed to confirm \Cref{con:chromLB} in full or at least for a `positive fraction' of all triangles.

Shkredov~\cite{shkredov2015some} used Fourier analysis to prove that every \textit{measurable} 2-coloring of the plane contains a monochromatic copy of each `not too close to equilateral' triangle. The points identifying normed triangles in this class would occupy roughly 80\% of the shaded area in \Cref{fig:allMono}. In a recent, not peer-reviewed manuscript~\cite{AI26}, OpenAI claims that it is enough to consider only \textit{measurable} 5-colorings to conclude that \textit{every} 5-coloring of the plane contains two points of the same color at distance 1. We wonder if these two ideas can be combined to confirm \Cref{con:chromLB} in full\footnote{It seems that a similar strategy was used in a complete solution to \Cref{con:chromLB} announced in a github repository~\cite{AI26_2} on the 19th of October 2026. The repository does not even contain a pdf-file with the proof, rather a sketch and a lean-verification code. Anton Bernshteyn also informed us in August 2026 that Hemanshu Kaul and he extended Shkredov's result and proved that every measurable 2-coloring of the plane contains a monochromatic copy of each non-equilateral triangle. Their manuscript, also containing a proof of \Cref{con:chromLB}, is currently in preparation.}.

\begin{figure}[htbp]
  \begin{subfigure}[b]{0.45\textwidth}
    \centering
      \resizebox{\linewidth}{!}{
      \begin{tikzpicture}
    \def\trisize{1}
    \trigrid{6}{7}
    \foreach \i in {-3,...,3}{
      \foreach \j in {-3,...,3}{
        \def\ci{6 * \i - 4 * \j}
        \def\cj{2 * \j}
        \node[circle, fill=red, minimum size=6pt, inner sep=0pt] at (\tricoord{\ci}{\cj+}) {};
        
        \node[diamond, fill=blue, minimum size=8pt, inner sep=0pt] at (\tricoord{\ci+1}{\cj}) {};
        \node[diamond, fill=blue, minimum size=8pt, inner sep=0pt] at (\tricoord{\ci+1}{\cj-1+}) {};
        \node[diamond, fill=blue, minimum size=8pt, inner sep=0pt] at (\tricoord{\ci}{\cj-1}) {};
        \node[diamond, fill=blue, minimum size=8pt, inner sep=0pt] at (\tricoord{\ci-1}{\cj}) {};
        \node[diamond, fill=blue, minimum size=8pt, inner sep=0pt] at (\tricoord{\ci-1}{\cj+1}) {};
        \node[diamond, fill=blue, minimum size=8pt, inner sep=0pt] at (\tricoord{\ci}{\cj+1}) {};

        \node[rectangle, fill=green, minimum size=6pt, inner sep=0pt] at (\tricoord{\ci+2}{\cj}) {};
        \node[rectangle, fill=green, minimum size=6pt, inner sep=0pt] at (\tricoord{\ci+2}{\cj-1}) {};
        \node[rectangle, fill=green, minimum size=6pt, inner sep=0pt] at (\tricoord{\ci+2}{\cj-2}) {};
        \node[rectangle, fill=green, minimum size=6pt, inner sep=0pt] at (\tricoord{\ci+1}{\cj-2}) {};
        \node[rectangle, fill=green, minimum size=6pt, inner sep=0pt] at (\tricoord{\ci}{\cj-2}) {};
        \node[rectangle, fill=green, minimum size=6pt, inner sep=0pt] at (\tricoord{\ci-1}{\cj-1}) {};
        \node[rectangle, fill=green, minimum size=6pt, inner sep=0pt] at (\tricoord{\ci-2}{\cj}) {};
        \node[rectangle, fill=green, minimum size=6pt, inner sep=0pt] at (\tricoord{\ci-2}{\cj+1}) {};
        \node[rectangle, fill=green, minimum size=6pt, inner sep=0pt] at (\tricoord{\ci-2}{\cj+2}) {};
        \node[rectangle, fill=green, minimum size=6pt, inner sep=0pt] at (\tricoord{\ci-1}{\cj+2}) {};
        \node[rectangle, fill=green, minimum size=6pt, inner sep=0pt] at (\tricoord{\ci}{\cj+2}) {};
        \node[rectangle, fill=green, minimum size=6pt, inner sep=0pt] at (\tricoord{\ci+1}{\cj+1}) {};

      }
    }
  \end{tikzpicture}}
    \captionsetup{justification=centering}
    \caption{Contains monochromatic and rainbow copies of $\Tr{1}{1}{2}$, but not of $\TrI{1}$ or $\TrI{2}$.}
  \end{subfigure}
  \hfill
  \begin{subfigure}[b]{0.45\textwidth}
    \centering
      \resizebox{\linewidth}{!}{
      \begin{tikzpicture}
    \def\trisize{1}
    \trigrid{6}{7}
    \def\xoff{0}
    \def\yoff{0}
    \foreach \i in {-4,...,4}{
      \foreach \j in {-4,...,4}{
        \def\cj{\j * 4 - \i}
        \def\ci{3 * \i}
        \node[circle, fill=red, minimum size=6pt, inner sep=0pt] at (\tricoord{\ci+\xoff}{\cj+\yoff}) {};
        \node[circle, fill=red, minimum size=6pt, inner sep=0pt] at (\tricoord{\ci+\xoff}{\cj+1+\yoff}) {};
        \node[circle, fill=red, minimum size=6pt, inner sep=0pt] at (\tricoord{\ci+1+\xoff}{\cj-1+\yoff}) {};
        \node[circle, fill=red, minimum size=6pt, inner sep=0pt] at (\tricoord{\ci+1+\xoff}{\cj+\yoff}) {};
      }
    }
    \def\xoff{1}
    \def\yoff{1}
    \foreach \i in {-4,...,4}{
      \foreach \j in {-4,...,4}{
        \def\cj{\j * 4 - \i}
        \def\ci{3 * \i}
        \node[diamond, fill=blue, minimum size=8pt, inner sep=0pt] at (\tricoord{\ci+\xoff}{\cj+\yoff}) {};
        \node[diamond, fill=blue, minimum size=8pt, inner sep=0pt] at (\tricoord{\ci+\xoff}{\cj+1+\yoff}) {};
        
        \node[diamond, fill=blue, minimum size=8pt, inner sep=0pt] at (\tricoord{\ci+1+\xoff}{\cj-1+\yoff}) {};
        \node[diamond, fill=blue, minimum size=8pt, inner sep=0pt] at (\tricoord{\ci+1+\xoff}{\cj+\yoff}) {};
      }
    }

    \def\xoff{-1}
    \def\yoff{3}
    \foreach \i in {-4,...,4}{
      \foreach \j in {-4,...,4}{
        \def\cj{\j * 4 - \i}
        \def\ci{3 * \i}
        \node[rectangle, fill=green, minimum size=6pt, inner sep=0pt] at (\tricoord{\ci+\xoff}{\cj+\yoff}) {};
        \node[rectangle, fill=green, minimum size=6pt, inner sep=0pt] at (\tricoord{\ci+\xoff}{\cj+1+\yoff}) {};
        \node[rectangle, fill=green, minimum size=6pt, inner sep=0pt] at (\tricoord{\ci+1+\xoff}{\cj-1+\yoff}) {};
        \node[rectangle, fill=green, minimum size=6pt, inner sep=0pt] at (\tricoord{\ci+1+\xoff}{\cj+\yoff}) {};
      }
    }
  \end{tikzpicture}}
    \captionsetup{justification=centering}
    \caption{Contains monochromatic and rainbow copies of $\TrI{1}$, but not of $\Tr{1}{1}{2}$.}
  \end{subfigure}
  \par
  \caption{Colorings of $\T$ showing that the na\"ive canonical analogue of \Cref{theorem:equi} fails.}
  \label{fig:equiCounter}
\end{figure}
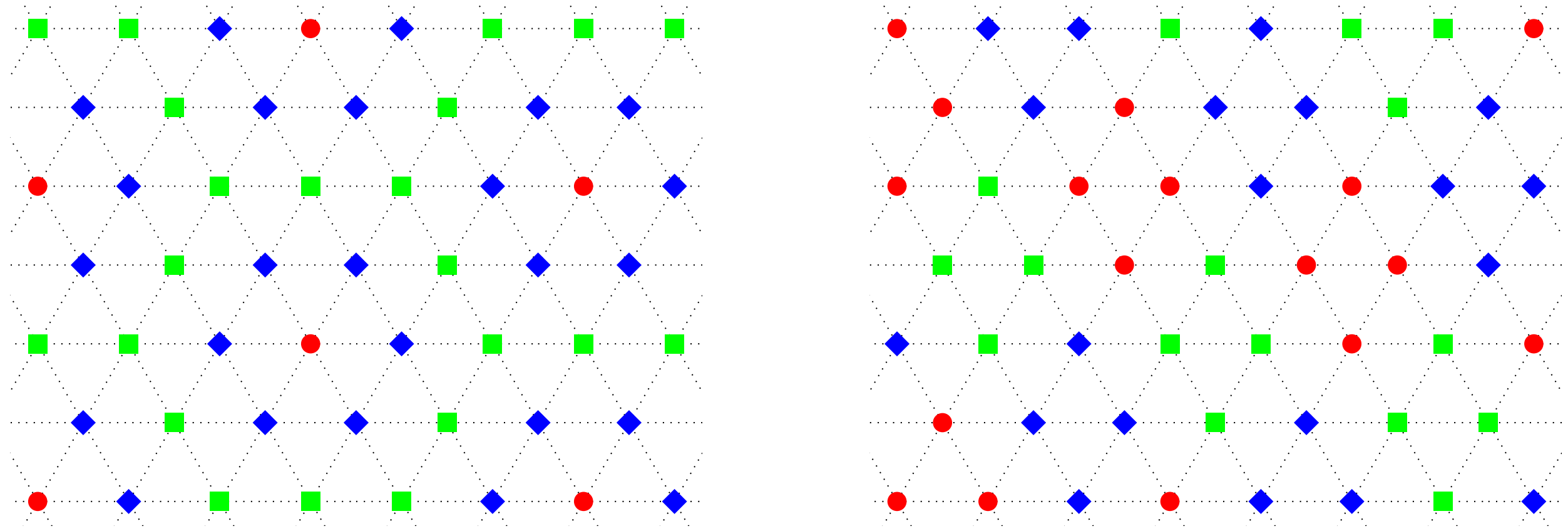

In this paper, we also confirmed \Cref{con:MR}, a canonical strengthening of \Cref{con:chromLB}, for the first 36 triangles. It would be interesting to extend this result at least to an infinite class of triangles, although it seems that different ideas are needed here too. Indeed, all known proofs of \Cref{con:chromLB} for infinite classes of triangles rely heavily on \Cref{theorem:equi}. We remark that all the auxiliary points in the proof of \Cref{theorem:equi} for a triangle with vertices on the triangular lattice $\T$ lie on $\T$ as well. The colorings of $\T$ in \Cref{fig:equiCounter} show that a canonical analogue of \Cref{theorem:equi}, even if exists, might be less elegant.

\vspace{2mm}
\noindent
\textbf{Acknowledgements.} We thank Anton Bernshteyn and Hemanshu Kaul for bringing the github repository~\cite{AI26_2} to our attention.

\appendix

\section{Supplementary data} \label{sec:appA}

\begin{longtable}{@{}lllll@{}}
\captionsetup{justification=centering}
\caption{Triangles $\TrO=\{\trip{0}{0},B,C\}$ satisfying $\R^2 \MRarrows \TrO$ \\ and thus $\R^2 \xrightarrow{2} \TrO$ found by the algorithm from Appendix~\ref{sec:app}.}
\label{tab:MR}\\
\toprule
$B$ & $C$ & Triangle & Shape & Belongs to classes \\
\midrule
\endfirsthead
\multicolumn{5}{@{}l}{\footnotesize\textit{Table continues}}\\
\toprule
$B$ & $C$ & Triangle & Shape & Belongs to classes \\
\midrule
\endhead
\midrule
\multicolumn{5}{r@{}}{\footnotesize\textit{Table continues on the next page}}\\
\endfoot
\bottomrule
\endlastfoot
$\trip{-3}{-1}$ & $\trip{1}{-3}$ & $\Tr{\sqrt{7}}{2 \sqrt{3}}{\sqrt{13}}$ & Acute & $-$ \\
$\trip{-3}{-2}$ & $\trip{-2}{2}$ & $\Tr{2}{\sqrt{19}}{\sqrt{21}}$ & Acute & $-$ \\
$\trip{-4}{-2}$ & $\trip{-3}{3}$ & $\Tr{3}{2 \sqrt{7}}{\sqrt{31}}$ & Acute & $-$ \\
$\trip{-4}{-3}$ & $\trip{2}{-4}$ & $\Tr{2 \sqrt{3}}{\sqrt{31}}{\sqrt{37}}$ & Acute & $-$ \\
$\trip{-3}{-1}$ & $\trip{-3}{0}$ & $\Tr{1}{3}{\sqrt{13}}$ & Obtuse & $-$ \\
$\trip{-3}{-1}$ & $\trip{1}{0}$ & $\Tr{1}{\sqrt{13}}{\sqrt{21}}$ & Obtuse & $-$ \\
$\trip{-3}{-1}$ & $\trip{3}{-3}$ & $\Tr{3}{\sqrt{13}}{2 \sqrt{7}}$ & Obtuse & $-$ \\
$\trip{-3}{-2}$ & $\trip{-1}{2}$ & $\Tr{\sqrt{3}}{\sqrt{19}}{2 \sqrt{7}}$ & Obtuse & $-$ \\
$\trip{-4}{1}$ & $\trip{2}{-4}$ & $\Tr{2 \sqrt{3}}{\sqrt{13}}{\sqrt{31}}$ & Obtuse & $-$ \\
$\trip{-4}{-1}$ & $\trip{2}{-2}$ & $\Tr{2}{\sqrt{21}}{\sqrt{31}}$ & Obtuse & $-$ \\
$\trip{-4}{-3}$ & $\trip{-3}{-2}$ & $\Tr{\sqrt{3}}{\sqrt{19}}{\sqrt{37}}$ & Obtuse & $-$ \\
$\trip{-4}{-3}$ & $\trip{-4}{-1}$ & $\Tr{2}{\sqrt{21}}{\sqrt{37}}$ & Obtuse & $-$ \\
$\trip{-4}{-2}$ & $\trip{1}{0}$ & $\Tr{1}{2 \sqrt{7}}{\sqrt{39}}$ & Obtuse & $-$ \\
$\trip{-4}{-3}$ & $\trip{-1}{2}$ & $\Tr{\sqrt{3}}{\sqrt{37}}{7}$ & Obtuse & $-$ \\
$\trip{-4}{-3}$ & $\trip{2}{-1}$ & $\Tr{\sqrt{3}}{\sqrt{37}}{2 \sqrt{13}}$ & Obtuse & $-$ \\
$\trip{-2}{0}$ & $\trip{-2}{1}$ & $\Tr{1}{\sqrt{3}}{2}$ & Right & (i), (iii), (ix), (xi), (xii) \\
$\trip{-2}{-1}$ & $\trip{-2}{1}$ & $\Tr{\sqrt{3}}{2}{\sqrt{7}}$ & Right & (ix) \\
$\trip{-2}{-2}$ & $\trip{-1}{1}$ & $\Tr{1}{2 \sqrt{3}}{\sqrt{13}}$ & Right & (ix) \\
$\trip{-2}{1}$ & $\trip{-1}{0}$ & $\Tr{1}{1}{\sqrt{3}}$ & Obtuse & (i), (v), (vii), (viii), (x), (xi) \\
$\trip{-2}{-1}$ & $\trip{-1}{-1}$ & $\Tr{1}{\sqrt{3}}{\sqrt{7}}$ & Obtuse & (ii), (viii), (xi) \\
$\trip{-2}{-1}$ & $\trip{1}{-1}$ & $\Tr{1}{\sqrt{7}}{3}$ & Obtuse & (viii) \\
$\trip{-2}{-2}$ & $\trip{-2}{-1}$ & $\Tr{1}{\sqrt{7}}{2 \sqrt{3}}$ & Obtuse & (i) \\
$\trip{-2}{0}$ & $\trip{1}{1}$ & $\Tr{\sqrt{3}}{2}{\sqrt{13}}$ & Obtuse & (ii), (x) \\
$\trip{-2}{-1}$ & $\trip{-1}{2}$ & $\Tr{\sqrt{3}}{\sqrt{7}}{\sqrt{13}}$ & Obtuse & (viii) \\
$\trip{-2}{-1}$ & $\trip{2}{-1}$ & $\Tr{\sqrt{3}}{\sqrt{7}}{4}$ & Obtuse & (i) \\
$\trip{-3}{-2}$ & $\trip{-3}{1}$ & $\Tr{\sqrt{7}}{3}{\sqrt{19}}$ & Obtuse & (x) \\
$\trip{-2}{-2}$ & $\trip{0}{1}$ & $\Tr{1}{2 \sqrt{3}}{\sqrt{19}}$ & Obtuse & (ii) \\
$\trip{-3}{-2}$ & $\trip{-2}{1}$ & $\Tr{\sqrt{3}}{\sqrt{13}}{\sqrt{19}}$ & Obtuse & (viii) \\
$\trip{-2}{-1}$ & $\trip{2}{0}$ & $\Tr{2}{\sqrt{7}}{\sqrt{21}}$ & Obtuse & (xi) \\
$\trip{-3}{-1}$ & $\trip{2}{-2}$ & $\Tr{2}{\sqrt{13}}{\sqrt{21}}$ & Obtuse & (viii) \\
$\trip{-4}{-1}$ & $\trip{-3}{-2}$ & $\Tr{1}{\sqrt{19}}{\sqrt{21}}$ & Obtuse & (viii) \\
$\trip{-3}{1}$ & $\trip{2}{2}$ & $\Tr{\sqrt{7}}{2 \sqrt{3}}{\sqrt{31}}$ & Obtuse & (viii) \\
$\trip{-4}{-3}$ & $\trip{1}{-2}$ & $\Tr{\sqrt{3}}{\sqrt{31}}{\sqrt{37}}$ & Obtuse & (viii) \\
$\trip{-4}{2}$ & $\trip{2}{3}$ & $\Tr{2 \sqrt{3}}{\sqrt{19}}{\sqrt{43}}$ & Obtuse & (viii) \\
$\trip{-4}{-4}$ & $\trip{-4}{-3}$ & $\Tr{1}{\sqrt{37}}{4 \sqrt{3}}$ & Obtuse & (i) \\
$\trip{-4}{-4}$ & $\trip{0}{1}$ & $\Tr{1}{4 \sqrt{3}}{\sqrt{61}}$ & Obtuse & (ii) \\
\end{longtable}

\begin{longtable}{@{}lllll@{}}
\captionsetup{justification=centering}
\caption{Triangles $\TrO=\{\trip{0}{0},B,C\}$ satisfying $\R^2 \xrightarrow{2} \TrO$ with \\ unresolved status of $\R^2 \MRarrows \TrO$ found by the algorithm from Appendix~\ref{sec:app}.}
\label{tab:MnoR}\\
\toprule
$B$ & $C$ & Triangle & Shape & Belongs to classes \\
\midrule
\endfirsthead
\multicolumn{5}{@{}l}{\footnotesize\textit{Table continues}}\\
\toprule
$B$ & $C$ & Triangle & Shape & Belongs to classes \\
\midrule
\endhead
\midrule
\multicolumn{5}{r@{}}{\footnotesize\textit{Table continues on the next page}}\\
\endfoot
\bottomrule
\endlastfoot
$\trip{-2}{-1}$ & $\trip{1}{-2}$ & $\Tr{\sqrt{3}}{\sqrt{7}}{\sqrt{7}}$ & Acute & $-$ \\
$\trip{-4}{-3}$ & $\trip{2}{-3}$ & $\Tr{\sqrt{7}}{6}{\sqrt{37}}$ & Acute & $-$ \\
$\trip{-4}{-4}$ & $\trip{-3}{2}$ & $\Tr{\sqrt{7}}{\sqrt{43}}{4 \sqrt{3}}$ & Acute & $-$ \\
$\trip{-2}{-1}$ & $\trip{1}{1}$ & $\Tr{\sqrt{3}}{\sqrt{7}}{\sqrt{19}}$ & Obtuse & $-$ \\
$\trip{-4}{-2}$ & $\trip{-4}{-1}$ & $\Tr{1}{\sqrt{21}}{2 \sqrt{7}}$ & Obtuse & $-$ \\
$\trip{-3}{-1}$ & $\trip{3}{-1}$ & $\Tr{\sqrt{7}}{\sqrt{13}}{6}$ & Obtuse & $-$ \\
$\trip{-3}{-2}$ & $\trip{3}{-2}$ & $\Tr{\sqrt{7}}{\sqrt{19}}{6}$ & Obtuse & $-$ \\
$\trip{-4}{-3}$ & $\trip{1}{-4}$ & $\Tr{\sqrt{13}}{\sqrt{21}}{\sqrt{37}}$ & Obtuse & $-$ \\
$\trip{-4}{-3}$ & $\trip{-2}{-4}$ & $\Tr{\sqrt{3}}{2 \sqrt{7}}{\sqrt{37}}$ & Obtuse & $-$ \\
$\trip{-3}{-2}$ & $\trip{-1}{3}$ & $\Tr{\sqrt{7}}{\sqrt{19}}{\sqrt{39}}$ & Obtuse & $-$ \\
$\trip{-4}{1}$ & $\trip{2}{2}$ & $\Tr{2 \sqrt{3}}{\sqrt{13}}{\sqrt{43}}$ & Obtuse & $-$ \\
$\trip{-4}{-4}$ & $\trip{-3}{-2}$ & $\Tr{\sqrt{7}}{\sqrt{19}}{4 \sqrt{3}}$ & Obtuse & $-$ \\
$\trip{-4}{-4}$ & $\trip{-3}{1}$ & $\Tr{\sqrt{7}}{\sqrt{31}}{4 \sqrt{3}}$ & Obtuse & $-$ \\
$\trip{-3}{-1}$ & $\trip{2}{2}$ & $\Tr{2 \sqrt{3}}{\sqrt{13}}{7}$ & Obtuse & $-$ \\
$\trip{-4}{-1}$ & $\trip{-1}{4}$ & $\Tr{\sqrt{13}}{\sqrt{21}}{7}$ & Obtuse & $-$ \\
$\trip{-4}{-2}$ & $\trip{1}{1}$ & $\Tr{\sqrt{3}}{2 \sqrt{7}}{7}$ & Obtuse & $-$ \\
$\trip{-4}{-3}$ & $\trip{-1}{3}$ & $\Tr{\sqrt{7}}{\sqrt{37}}{3 \sqrt{7}}$ & Obtuse & $-$ \\
$\trip{-2}{-1}$ & $\trip{-2}{2}$ & $\Tr{2}{\sqrt{7}}{3}$ & Acute & (x) \\
$\trip{-3}{-2}$ & $\trip{2}{-3}$ & $\Tr{\sqrt{7}}{\sqrt{19}}{\sqrt{21}}$ & Acute & (xi) \\
$\trip{-3}{-3}$ & $\trip{-3}{2}$ & $\Tr{\sqrt{7}}{5}{3 \sqrt{3}}$ & Acute & (i) \\
$\trip{-4}{-2}$ & $\trip{3}{-4}$ & $\Tr{\sqrt{13}}{2 \sqrt{7}}{\sqrt{39}}$ & Acute & (xi) \\
$\trip{-4}{-4}$ & $\trip{-4}{3}$ & $\Tr{\sqrt{13}}{4 \sqrt{3}}{7}$ & Acute & (i) \\
$\trip{-3}{-2}$ & $\trip{1}{-2}$ & $\Tr{\sqrt{3}}{4}{\sqrt{19}}$ & Right & (ix) \\
$\trip{-4}{-2}$ & $\trip{1}{-2}$ & $\Tr{\sqrt{3}}{5}{2 \sqrt{7}}$ & Right & (ix) \\
$\trip{-3}{-3}$ & $\trip{-1}{1}$ & $\Tr{1}{3 \sqrt{3}}{2 \sqrt{7}}$ & Right & (ix) \\
$\trip{-3}{-3}$ & $\trip{-2}{2}$ & $\Tr{2}{3 \sqrt{3}}{\sqrt{31}}$ & Right & (ix) \\
$\trip{-4}{-3}$ & $\trip{-4}{2}$ & $\Tr{2 \sqrt{3}}{5}{\sqrt{37}}$ & Right & (ix) \\
$\trip{-4}{4}$ & $\trip{-3}{-3}$ & $\Tr{4}{3 \sqrt{3}}{\sqrt{43}}$ & Right & (ix) \\
$\trip{-4}{-4}$ & $\trip{-1}{1}$ & $\Tr{1}{4 \sqrt{3}}{7}$ & Right & (ix) \\
$\trip{-3}{-1}$ & $\trip{-2}{3}$ & $\Tr{\sqrt{7}}{\sqrt{13}}{\sqrt{21}}$ & Obtuse & (xi) \\
$\trip{-3}{2}$ & $\trip{2}{2}$ & $\Tr{\sqrt{7}}{2 \sqrt{3}}{5}$ & Obtuse & (i) \\
$\trip{-3}{-1}$ & $\trip{2}{-1}$ & $\Tr{\sqrt{3}}{\sqrt{13}}{5}$ & Obtuse & (i) \\
$\trip{-3}{-3}$ & $\trip{-3}{-1}$ & $\Tr{2}{\sqrt{13}}{3 \sqrt{3}}$ & Obtuse & (i) \\
$\trip{-3}{-3}$ & $\trip{-3}{1}$ & $\Tr{\sqrt{7}}{4}{3 \sqrt{3}}$ & Obtuse & (i) \\
$\trip{-3}{-3}$ & $\trip{-3}{-2}$ & $\Tr{1}{\sqrt{19}}{3 \sqrt{3}}$ & Obtuse & (i) \\
$\trip{-3}{0}$ & $\trip{1}{2}$ & $\Tr{\sqrt{7}}{3}{2 \sqrt{7}}$ & Obtuse & (xii) \\
$\trip{-4}{-3}$ & $\trip{-3}{1}$ & $\Tr{\sqrt{7}}{\sqrt{21}}{\sqrt{37}}$ & Obtuse & (xi) \\
$\trip{-3}{-3}$ & $\trip{0}{1}$ & $\Tr{1}{3 \sqrt{3}}{\sqrt{37}}$ & Obtuse & (ii) \\
$\trip{-3}{-1}$ & $\trip{2}{1}$ & $\Tr{\sqrt{7}}{\sqrt{13}}{\sqrt{39}}$ & Obtuse & (xi) \\
$\trip{-4}{0}$ & $\trip{1}{2}$ & $\Tr{\sqrt{7}}{4}{\sqrt{39}}$ & Obtuse & (viii) \\
$\trip{-4}{-2}$ & $\trip{-2}{3}$ & $\Tr{\sqrt{7}}{2 \sqrt{7}}{\sqrt{39}}$ & Obtuse & (xii) \\
$\trip{-4}{-1}$ & $\trip{3}{-2}$ & $\Tr{\sqrt{7}}{\sqrt{21}}{\sqrt{43}}$ & Obtuse & (xi) \\
$\trip{-4}{-4}$ & $\trip{-4}{1}$ & $\Tr{\sqrt{13}}{5}{4 \sqrt{3}}$ & Obtuse & (i) \\
$\trip{-4}{2}$ & $\trip{3}{2}$ & $\Tr{2 \sqrt{3}}{\sqrt{19}}{7}$ & Obtuse & (i) \\
$\trip{-3}{-3}$ & $\trip{0}{2}$ & $\Tr{2}{3 \sqrt{3}}{7}$ & Obtuse & (ii) \\
$\trip{-4}{-1}$ & $\trip{2}{1}$ & $\Tr{\sqrt{7}}{\sqrt{21}}{2 \sqrt{13}}$ & Obtuse & (xi) \\
$\trip{-4}{-1}$ & $\trip{4}{-3}$ & $\Tr{\sqrt{13}}{\sqrt{21}}{2 \sqrt{13}}$ & Obtuse & (xii) \\
$\trip{-4}{1}$ & $\trip{3}{2}$ & $\Tr{\sqrt{13}}{\sqrt{19}}{\sqrt{57}}$ & Obtuse & (xi) \\
$\trip{-4}{-2}$ & $\trip{3}{-1}$ & $\Tr{\sqrt{7}}{2 \sqrt{7}}{\sqrt{57}}$ & Obtuse & (xii) \\
$\trip{-3}{-3}$ & $\trip{1}{2}$ & $\Tr{\sqrt{7}}{3 \sqrt{3}}{\sqrt{61}}$ & Obtuse & (viii) \\
$\trip{-4}{0}$ & $\trip{3}{3}$ & $\Tr{4}{3 \sqrt{3}}{\sqrt{79}}$ & Obtuse & (ii) \\
$\trip{-2}{0}$ & $\trip{1}{0}$ & $\Tr{1}{2}{3}$ & Degenerate & (iii), (viii), (xii) \\
\end{longtable}

\section{Program code} \label{sec:app}

{\captionsetup[lstlisting]{justification=centering}
\begin{lstlisting}[language=Python,
  captionpos=t,
  title={Python code proving \protect\Cref{th:2}. To prove \protect\Cref{th:3} instead, \\ remove the restriction \texttt{c < 2} in Line 77.}
  %label={lst:algoMono}
]
#!/usr/bin/env python3
from operator import itemgetter
from itertools import product, combinations
from z3 import Int, Solver, And, Not, sat

# Operation for adding tuples component-wise.
def t_add(tuple_1, tuple_2):
    return tuple(map(sum, zip(tuple_1, tuple_2)))

# Offset a list of points by some fixed tuple.
def offset_points(points, offset):
    return [t_add(point, offset) for point in points]

# Appends a candidate to main_list
# if candidate coordinates are within some radius
def try_append(main_list, candidate, radius):
    for coord in candidate:
        if abs(coord[0]) > radius or abs(coord[1]) > radius:
            return
    main_list.append(candidate)

# Returns all points of the radius-box T_{radius}.
def build_triangular_lattice(radius):
    return product(range(-radius, radius+1), repeat=2)

# Rotates a point coord around the origin by 60 * angle degrees.
def rotate_by(coord, angle=1):
    if angle == 0:
        return coord
    else:
        rotated_by_1 = (-coord[1], coord[0] + coord[1]) 
        return rotate_by(rotated_by_1, angle=angle-1)

# Rotates a set of points around the origin by 60 * angle degrees.
def rotate_points(points, angle=1):
    if angle == 0:
        return points
    else:
        rotated_by_1 = [rotate_by(coord, angle=1) for coord in points]
        return rotate_points(rotated_by_1, angle=angle-1)

# Returns a list of all congruent copies of some
# triangle (0, point_1, point_2) within T_{radius}.
def find_all_triangles(radius, point_1, point_2):
    triangles = []
    pts = build_triangular_lattice(radius)
    base_triangle = [(0,0), point_1, point_2]
    flipped_triangle = [(0,0), (point_1[1], point_1[0]), (point_2[1], point_2[0])]
    for pt in pts:
        for rotation in range(6):
            try_append(triangles, offset_points(rotate_points(base_triangle, angle=rotation), pt), radius)
            try_append(triangles, offset_points(rotate_points(flipped_triangle, angle=rotation), pt), radius)

    # Remove duplicates:
    triangles = [tuple(sorted(tri, key=itemgetter(0,1))) for tri in triangles]
    triangles = list(set(triangles))
    return triangles

# Returns a unique integer for each coord in T_{radius}.
def vtx_index(coord, radius):
    shifted_coord = t_add(coord, (radius, radius))
    return shifted_coord[0] * ( 2 * radius + 1) + shifted_coord[1]

# Returns the unique coord determined by its vtx_index.
def rev_vtx_index(idx, radius):
    y = idx % (2 * radius + 1)
    x = idx // (2 * radius + 1)
    return t_add((x,y), (-radius, -radius))

# Returns True when some coloring of T_{radius} is found that
# avoids every monochromatic configuration in all_triangles.
def solve_sat(all_triangles, radius, diff_pairs=None):
    n = (2*radius + 1)**2  # number of vertices
    color_vars = [Int(f"c_{i}") for i in range(n)]
    solver = Solver()
    for c in color_vars:
        solver.add(c >= 0, c < 2)
    # enforce "different color" constraints
    if diff_pairs:
        for pair in diff_pairs:
            i = vtx_index(pair[0], radius)
            j = vtx_index(pair[1], radius)
            solver.add(color_vars[i] != color_vars[j])
    # Every congruent copy must not be monochromatic
    for tri in all_triangles:
        idxs = [vtx_index(v, radius) for v in tri]
        c0, c1, c2 = color_vars[idxs[0]], color_vars[idxs[1]], color_vars[idxs[2]]
        # Forbid monochromatic and rainbow copies of tri
        solver.add(Not(Or(And(c0 != c1, c0 != c2, c1 != c2), And(c0 == c1, c1==c2, c2==c0))))
    # wLog: Fix first vertex to 0 to break symmetry
    solver.add(color_vars[0] == 0)
    # Solve
    if solver.check() != sat:
        # No coloring exists
        return False
    # Some coloring exists
    else:
        return True

r=6
max_dist = 4
# Check all triangles within T_{max_dist}
for idx, (p1,p2) in enumerate(combinations(build_triangular_lattice(max_dist), 2)):
    print(f"--- ({idx}) Now checking {p1}, {p2} ---")
    triangles = find_all_triangles(r,p1,p2)
    for i in range(-r+1, r+1):
        solve_result = solve_sat(triangles, radius=r, diff_pairs=[[(-r,0), (i,0)]])
        if not solve_result:
            print(f"Triangle {p1}, {p2} has no coloring with distance {abs(r + i)}.")
            break
\end{lstlisting}}

\end{document}